\documentclass{amsart}

\usepackage{macros}
\standardmargins\standardpackages\standardrenews\standardlabeling\standardmakes
\colorcommentstrue

\draftfalse

\authormumtaz
\authorsimon
\authordavid

\DeclareMathOperator{\mat}{Mat}
\newcommand{\rk}{\mathrm{rank}}
\newcommand{\norm}[1]{\ensuremath{\left\Vert #1 \right\Vert}}
\newcommand{\abs}[1]{\ensuremath{\left\vert #1 \right\vert}}

\newcommand{\wtilde}{\widetilde}
\renewcommand{\subset}{\subseteq}

\newcommand{\Badmn}{ {\rm Bad}(m,n)}
\newcommand{\Bu}{ {\rm Bad}(m,n,u)}
\newcommand{\Balpha}{ {\rm Bad}_\bfalpha(m,n)}

\newcommand{\Bbo}{ {\rm Bad}_\0(m,n)}

\begin{document}
\title{Metrical theorems on systems of affine forms}

\begin{abstract}

In this paper we discuss metric theory associated with the affine (inhomogeneous) linear forms in the so called doubly metric settings within the classical and the mixed setups. We consider the system of affine forms given by $\qq\mapsto \qq X+\bfalpha$, where $\qq\in\Z^m$ (viewed as a row vector), $X$ is an $m\times n$ real matrix and $\bfalpha\in \R^n$. The classical setting refers to the  ${\rm dist}(\qq X+\bfalpha, \Z^m)$ to measure the closeness of the integer values of the system $(X, \bfalpha)$ to integers. The absolute value setting is obtained by replacing ${\rm dist}(\qq X+\bfalpha, \Z^m)$ with ${\rm dist}(\qq X+\bfalpha, \0)$; and the more general mixed settings are obtained by replacing ${\rm dist}(\qq X+\bfalpha, \Z^m)$ with ${\rm dist}(\qq X+\bfalpha, \Lambda)$, where $\Lambda$ is a subgroup of $\Z^m$.

We prove the  Khintchine--Groshev  and Jarn\'ik type theorems for the mixed affine forms and Jarn\'ik type theorem for the classical affine forms. We further prove that the sets of badly approximable affine forms,  in both the classical and mixed settings,  are hyperplane winning.  The latter result,  for the classical setting, answers a question raised by Kleinbock (1999).

%prove the Khintchine--Groshev  and Jarn\'ik type theorems for systems of linear forms in the so called doubly metric settings. Further,  we establish complete Khintchine--Groshev and Jarn\'ik--Schmidt type theorems for inhomogeneous small linear forms (i.e. ``absolute value approximation'') in the so-called doubly metric case, in which the inhomogeneous parameter is not fixed. In particular, we prove that the appropriate analogue of the set of badly approximable affine forms is hyperplane winning. The analogue of the Khintchine--Groshev theorem takes one of two different forms depending on a certain inequality between the dimensions appearing in the setup.
\end{abstract}

\maketitle

\section{Background and statements of results}

\subsection{Classical affine forms}

Let $\psi:\N\to\R^+ $ be a function tending to $0$ at infinity, referred to as an \emph{approximation function}. Let $W(m,n;\psi)$ denote the set of pairs $( X, \bfalpha) \in \R^{mn}\times \R^n$ for which the system of inequalities
\begin{equation}\label{eq1}
|q_1 x_{1i}+q_{2}x_{2i}+\dots+q_{m}x_{mi}-\alpha_i - p_i| \leq \psi(|\qq|)\text{\quad{}for\quad}1\leq i\leq n
\end{equation}
is satisfied for infinitely many $\qq\in \Z ^{m}\setminus\{\0\}$  and $\pp\in \Z^n$. On the right-hand side, $|\qq|$ denotes the max norm of $\qq$. In what follows, the system $$ \{q_1x_{1i}+ \dots+ q_{m}x_{mi} \; : \; i = 1,\ldots,n\}$$ of $n$ linear forms in the $m$ variables $q_1,\ldots,q_m$ will be written more concisely as $\qq X$, where the matrix $X$ is regarded as a point in $ \R^{mn}$. If the vector $\bfalpha$ is fixed then the study of the measure and dimension of the fiber $$W_
\bfalpha \df \{X \in \R^{mn} : (X,\bfalpha)\in W(m,n;\psi)\}$$ is referred as the \emph{singly metric} theory, while the study of the entire set $W(m, n;\psi)$ is called the \emph{doubly metric} theory. Both the singly metric and the doubly metric theories are considered part of \emph{inhomogeneous Diophantine approximation}. Inhomogeneous Diophantine approximation is in some respects different from homogeneous Diophantine approximation, which corresponds to the case $\bfalpha=\0$. In particular,  in the doubly metric case, the results are sometimes sharper and easier to prove than singly metric case due to the extra variable involved which offers an extra degree of freedom. 

The most fundamental result in metric Diophantine approximation is the Khintchine--Groshev theorem which gives an elegant answer to the question of the  `size' of the set in terms of Lebesgue or Hausdorff measure.  For the modernised version of the Khintchine--Groshev theorem for $W_\0(m,n;\psi)$ we refer to \cite{BDV} and for $W_\bfalpha(m,n;\psi)$ we refer to \cite{AllenBeresnevich}.  The doubly metric version of the Khintchine--Groshev theorem not only requires weaker assumptions than the homogeneous analogue  but is also considerably easier to prove, see \cite[Theorem 15]{Sprindzuk}, or \cite[Theorem VII.II]{Cassels} for the special case of simultaneous approximation.

Another set of significant interest is in a sense complementary to $W(m,n;\psi)$, namely the set of inhomogeneous badly approximable affine forms. First we define the classical set of badly approximable affine forms as

\[{\rm Bad}(m,n)=\{( X, \bfalpha) \in
\R^{mn}\times \R^n: \inf_{\qq\in\Z^m\setminus \{\0\}} |\qq|^{m/n}\|\qq X-\bfalpha\|>0\}
\]
and for each $\bfalpha\in\R^n$ we let
\[
\Balpha = \{X\in \R^{mn} : (X,\bfalpha)\in {\rm Bad}(m,n)\}.
\]

In his seminal paper \cite{Kleinbock4}, Kleinbock used ideas and techniques from the theory of dynamical systems to prove a doubly metric  result, namely that the set ${\rm Bad}(m,n)$ is of full Hausdorff dimension. Essentially, his method is based on a deep connection between badly approximable systems of linear forms and orbits of certain lattices in Euclidean space under an appropriate action. One possible strengthening of such dimension results is to show that the sets in question are winning for what is called Schmidt's game. Schmidt's game was introduced by Schmidt in \cite {Schmidt1} and later used by himself in \cite{Schmidt2} to prove that the set ${\rm Bad}_\0(m,n)$ of homogeneous badly approximable matrices is of full dimension. Actually this possible strengthening in the setting of this paper was conjectured by Kleinbock (see \cite[pp. 101]{Kleinbock4})\footnote{In the quoted text, we have adopted notation of our paper to avoid ambiguity. We would like to thank Lovy Singhal for bringing this conjecture to our attention.} as follows:

``\emph{It seems natural to conjecture that ${\rm Bad}(m,n)$ is  a wining subset of $\R^{mn}\times \R^n$,  and,  moreover,  that $\Balpha $ is  a  winning  subset  of  $\R^{mn}$ for  every $\bfalpha\in\R^n$.  This  seems  to  be  an interesting and challenging  problem in metric number theory.}''

Einsiedler and Tseng in \cite{EinsiedlerTseng} proved the latter half of the Kleinbock's conjecture by proving that each fiber $\Balpha$ is winning for Schmidt's game, which in-return is a considerable strengthening of just proving that the set in question has full Hausdorff dimension.  We prove the former half of the conjecture in a stronger form.

There have been a lot of developments recently not only in the usage of Schmidt's game but also progress has been made in introducing variants of this game. The \emph{hyperplane game} is a variant of a game defined by McMullen \cite{McMullen} which in turn is a variant of Schmidt's game \cite{Schmidt1}. It was originally defined by Broderick, Fishman, Kleinbock, Reich, and Weiss \cite{BFKRW}. Sets which are winning for the hyperplane game are called \emph{hyperplane winning}. In \cite{BFS1}, Broderick, Fishman, and the third-named author proved that each fiber $\Balpha$ is hyperplane winning, which is a strengthening of the result of Einsiedler and Tseng mentioned above. 

\smallskip

We prove the following result and thus fully settle Kleinbock's conjecture.

\begin{theorem}
\label{thm:BadHAW}
The set ${\rm Bad}(m,n)$ is hyperplane winning.
\end{theorem}

Note that the hyperplane winning property is valid even if the games are played on certain fractals. Hence this theorem also implies that the set ${\rm Bad}(m,n)$ intersection with a large class of fractal sets has full Hausdorff dimension, that is, the dimension of the fractal sets. We refer the reader to \cite{BFKRW} for details. 

%\commumtaz{ maybe we state what we actually mean by taking intersection with the fractals}

%A natural question which presents itself at this stage is the nature of the the metric structure of $W(m, n;\psi)$ if distance to the nearest integer in \eqref{eq1}, $\| \, \cdot \, \|$ is replaced by absolute value norm $| \,\cdot \, |$, \emph{i.e.}, by replacing the nearest integer vector in \eqref{eq1} to with the zero vector and thus considering system of linear forms which are simultaneously close to the origin. The question may be asked in the broader context if the set consists of systems of linear forms in which the
% forms of approximation are mixed and some forms are required to be small and others close to an integer.

\subsection{Mixed affine forms}

An extension of the classical problems which has drawn some attention recently is the following: what happens if one requires the pair $(\pp,\qq)$ appearing in the definition of the set $W(m, n;\psi)$ to belong to some fixed set of integer vectors? An example of such a result is a theorem of Dani, Laurent, and Nogueira \cite{DaniLauNog}, which is a doubly metric analogue of the Khintchine--Groshev theorem for the set $W(m, n;\psi)$ under the additional restriction that, with respect to the direct sum decomposition of $\R^n$ corresponding to a fixed and sufficiently coarse partition $\pi$ of the set of coordinates $\{1,\ldots,m+n\}$, the components of the integer vector $(\pp,\qq)$ are all primitive. In this paper we consider another natural  restriction, namely that the vector $\pp$ must lie in some sublattice of $\Z^n$.

More precisely, fix an integer $0 \leq u \leq n$, and consider sublattice  $\Z^u\times \{\0\}^{n-u}\subset \Z^n$. Let $W(m,n,u;\psi)$ be the set of pairs $(X, \bfalpha)\in \R^{mn}\times \R^n$ such that
\begin{multline}\label{mixeddef}
\max\{|\qq\cdot
\xx^{(1)}-\alpha_1-p_1|,\ldots,|\qq\cdot \xx^{(u)}-\alpha_u-p_u|, \\ |\qq\cdot
\xx^{(u+1)}-\alpha_{u+1}|, \ldots,|\qq\cdot
\xx^{(n)}-\alpha_n|\} \leq \psi(|\qq|)
\end{multline}
\noindent is satisfied for infinitely many $\qq\in\Z^m\setminus \{\0\}$ and $\pp\in \Z^u\times \{\0\}^{n-u}$. Here, $\xx^{(1)},\ldots,\xx^{(n)}$ are the column vectors of $X$.  These linear forms are referred to as the  \emph{absolute affine linear forms}.  The metrical theory for these sets and these particular sublattices contains the results for general sublattices, as the analogous set for a general sublattice is easily seen to be bi-Lipschitz equivalent to a set which can be bounded from inside and outside by sets of the form $W(m,n,u;\widetilde\psi)$ for appropriate functions $\widetilde\psi$ (see e.g. \cite{HussainKristensen}). Note that
\[
W(m,n;\psi) = W(m,n,n;\psi).
\]
% \comdavid{There was an incorrect statement here claiming bi-Lipschitz equivalence. However, applying a bi-Lipschitz transformation changes both the set of integers and the height function, so it is necessary to change the function $\psi$ to compensate. If $\psi$ is doubling then it suffices to multiply $\psi$ by a constant to compensate (large for the upper boound and small for the lower bound), but in general the relation between $\psi$ and $\widetilde\psi$ is slightly annoying to write down.}

Our new results on the metric properties of the sets $W(m,n,u;\psi)$ and the fibers
\[
W_\bfalpha(m,n,u;\psi) = \{X\in \R^{mn}: (X,\bfalpha)\in W(m,n,u;\psi)\}
\]
are proved using the methods developed in \cite{HussainKristensen2}, including the appeal to the important tools of Beresnevich and Velani, the mass transference principle \cite{BeresnevichVelani} and its generalization to linear forms \cite{AllenBeresnevich,BeresnevichVelani4}. The adaptation of the methods to the present setup requires some work, so rather than stating just the differences with the manuscript \cite{HussainKristensen2}, we have chosen to present the present work in a more self-contained manner.

In the case $u=0$, $\bfalpha= \0$, the set 
\begin{equation*}W_\0(m,n,0;\psi):=\left\{X\in\R^{mn}: \begin{array}{l} \max\{|\qq\cdot
\xx^{(1)}|, \ldots,|\qq\cdot
\xx^{(n)}|\}
\leq \psi(|\qq|) \\  {\rm for \ infinitely \ many} \ \qq \in \Z^{m}\setminus\{\0\} \end{array} \right\}
\end{equation*}
 is well studied by various authors in \cite{DickinsonHussain, Dickinson, FHKL, HussainKristensen2, HussainLevesley}. The set arose in the literature in part due to its connections with the Kolmogorov--Arnold--Moser theory \cite{Arnold, DRV_KAM}, linearization of germs of complex analytic diffeomorphisms of $\C^m$ near a fixed point \cite{Arnold, DRV4}, operator theory \cite{DGY1, DGY2} and recently discovered applications in signal processing \cite{MHMK, OrdUri}.

Note that $W(1,n,0;\psi)= \{(\0, \0)\}.$ Indeed, any $(x_1, \ldots, x_n, \alpha_1, \ldots, \alpha_n)\in W(1,n,0;\psi)$ must satisfy the inequality
\[
\underset{1\leq i\leq n}{\max}\left|x_i-\frac{\alpha_i}{q}\right|<\frac{\psi(q)}{q}
\]
for infinitely many $q\in \Z\setminus \{0\}$. Since $\alpha_i/q\to 0$ and $\psi(q)/q\to 0$ as $q\to\infty$, taking the limit as $q\to\infty$ yields $x_1 = \ldots = x_n = 0$, and since $\psi(q)\to 0$ as $q\to\infty$, multiplying by $q$ and then taking the limit yields $\alpha_1 = \ldots = \alpha_n = 0$. So $(\xx,\bfalpha)\in W(1,n,0;\psi)$ only when $(\xx,\bfalpha) = (\0,\0)$.

In a similar way, it can be readily verified that $W(1,n,u;\psi)= W(1,u;\psi)\times \{(\0,\0)\}$ for $u<n$ and thus
\[
\dim W(1,n,u; \psi) \ = \ \dim W(1, u; \psi).
\]
Consequently, throughout this paper we assume that $m>1$.

%\comsimon{We need a passage here on badly approximable systems in both settings.}\commumtaz{I have reshuffled the sections, probably we don't need the passage here any more}

\medskip

\noindent{\bf Notation.} Before proceeding with the description of our main results, we introduce some notation. The Vinogradov symbols $\lessless$ and $\gtrgtr$ will be used to indicate an inequality with an unspecified positive multiplicative constant depending only on $m$ and $n$. If $a\lessless  b$ and $a\gtrgtr b$ we write $a\asymp b$, and say that the quantities $a$ and $b$ are comparable. A \emph{dimension function} is an increasing continuous function $f:\R^+\rightarrow \R^+$ such that $f(r)\to 0$ as $r\to 0$. Throughout the paper, $\HH^f$ denotes the $f$--dimensional Hausdorff measure which will be fully defined in section \ref{hm}. Finally, for convenience, if $\psi$ is a given approximation function, then we write
$$\Psi(q):=\frac{\psi(q)}{q}.$$
The Hausdorff dimension of a set $E$ will be denoted by $\dim E$.

\medskip

%\subsection{Results}
The measure-theoretic results below for absolute value approximation crucially depend upon whether $m+u>n$ or $m+u\leq n$. The reason for a dichotomy should be clear from Lemma \ref{dimlemma} below.
\begin{theorem}\label{thm2}
Let $m+u>n$ and let $\psi$ be an approximating function. Let $f$ be a dimension function such that the maps $r \mapsto r^{-(m+1)n}f(r)$ and $r \mapsto r^{-(m+u+1-n)n}f(r)$ are both monotonic. Then for every nonempty open set $U \subset \R^{mn}\times \R^n$ we have
\[
\HH^{f}\big(W(m,n,u;\psi)\cap U\big)
= \left\{
\begin{array}{cl}
0 &
{\rm \ if} \ \sum
\limits_{q=1}^{\infty}f(\Psi (q))\Psi(q)^{-mn}q^{m+u-1}<\infty,\\
&\\
\HH^f(U)& {\rm \ if} \ \sum
\limits_{q=1}^{\infty}f(\Psi (q))\Psi(q)^{-mn}q^{m+u-1}=\infty.
\end{array}\right.
\]
\end{theorem}

Note that we have not assumed the approximating function $\psi$ to be monotonic for any $m$ and $n$. This is possible due to the extra degree of freedom offered by the variable $\bfalpha$. This is a doubly metric inhomogeneous version of the set considered in \cite{HussainKristensen2, HussainLevesley}, where an analogous result was obtained for the homogeneous fiber obtained by fixing $\bfalpha=\0$.

As in most of the statements the convergence part is reasonably straightforward to establish and is free from any assumptions on $m, u, n$, the approximation function, and the dimension function $f$. It is the divergence statement which constitutes the main substance and this is where conditions come into play.

The requirement that $r\mapsto r^{-(m+1)n}f(r)$ and other functions must be monotonic is a natural and not particularly restrictive condition. Essentially, the condition ensures that the Hausdorff measures cannot be too degenerate compared with the Lebesgue measure of the ambient Euclidean space, which is of dimension $(m+1)n$. In the case where $f(r) := r^{(m+1)n}$ the Hausdorff measure $\HH^f$ is proportional to the standard $(m+1)n$--dimensional Lebesgue measure of $\R^{mn}\times \R^n$ and the resulting special case of Theorem \ref{thm2} is the natural analogue of the Khintchine--Groshev theorem for $W(m,n,u;\psi)$:

\begin{corollary}\label{cor2}
Let $m+u>n$ and let $\psi$ be an approximating function. Then for every nonempty open set $U \subset \R^{mn}\times \R^n$,
\[
\big| W(m,n,u;\psi) \cap U \big|_{(m+1)n} = \left\{
\begin{array}{cl}
0 &
{\rm \ if} \
\sum\limits_{q=1}^\infty{}\Psi(q)^n q^{m+u-1}<\infty,\\
&\\
|U|_{(m+1)n}& {\rm \ if} \
\sum\limits_{q=1}^\infty{}\Psi(q)^n q^{m+u-1}=\infty.
\end{array}\right.\]
\end{corollary}

\medskip

From this corollary, it can easily be seen that for $\tau>\frac{m+u}{n}-1$, the set $W(m,n,u;\tau) \df W(m,n,u; q\mapsto q^{-\tau})$ is a null set.

The following corollary of Theorem \ref{thm2} gives the Hausdorff dimension of the set $W(m,n,u;\psi)$ as well as its measure with respect to the dimension functions $f(r)=r^s$, $s>0$.
\begin{corollary}\label{co1}
Let $m+u>n$, and let $\psi$ be an approximating function. Then for all $s\geq 0$  and for every nonempty open set $U \subset \R^{mn}\times\R^n$,
\[
\HH^s\big(W(m,n,u;\psi) \cap U\big) \ = \left\{
\begin{array}{cl}
0& {\rm \ if} \ \
\sum\limits_{q=1}^{\infty}\Psi(q)^{s-mn} q^{m+u-1}<\infty, \, \\
&\\
\HH^s(U) & { \rm \ if} \ \
\sum\limits_{q=1}^{\infty} \Psi(q)^{s-mn} q^{m+u-1}=\infty.
\end{array}
\right.
\]
Consequently, if $U$ is nonempty then
\[
\dim\big(W(m,n,u;\psi)\cap U\big) = \sup \left\{mn \leq s \leq (m+1)n : \sum_{q=1}^{\infty} \Psi(q)^{s-mn}q^{m+u-1} = \infty\right\}.
\]
\end{corollary}

Finally, for completeness, the dimension result for $W(m,n,u; \tau)$ is given for $m+u>n$. This follows directly from the previous corollary.

\begin{corollary} \label{cor3}
For $m+u>n$ and every nonempty open set $U \subset \R^{mn}\times \R^n$,
\[
\dim\big(W(m,n,u;\tau) \cap U\big) =\left\{\begin{array}{ll}
mn+\frac{m+u}{\tau+1} & {\rm \ if} \ \tau\geq\frac{m+u}{n}-1,\\
&\\
(m+1)n & {\rm \ if}\ \tau\le\frac{m+u}{n}-1.
\end{array}\right.\]
\end{corollary}

In fact, Corollary \ref{co1} gives more than simply the Hausdorff dimension. It also give the Hausdorff measure at the critical exponent which is infinity in this case, except when $W(m,n,u,\tau)$ is of full dimension.

\medskip

To state an analogue of Theorem \ref{thm2} for the case $m+u\leq n$, we will need to modify the conditions on the dimension function $f$. This change is due to the fact that if $m+u\leq n$, then the set $W(m,n,u;\psi)$ is contained in an algebraic variety of dimension strictly lower than $(m+1)n$. 
%\comdavid{Earlier there was written something about an over-determined system which I did not understand. Normally an over-determined system of equations is one that has more equations than variables, and I did not see anything like that here. So I removed the terminology about over-determined systems.}

To see this, first consider the case $m=n, u=0$ and fix $(X,\bfalpha) \in W(m,m,0;\psi)$ such that $\det X\neq 0$. Multiplying the defining inequalities (\ref{mixeddef}) by $X^{-1}$ shows that we must have
\[|\qq-\bfalpha X^{-1}|\leq C(X)\psi(|\qq|)\]
for infinitely many $\qq$. As $|\qq|$ tends to infinity, the left-hand side tends to infinity while the right-hand side tends to zero, and we get a contradiction. Hence, any element $(X,\bfalpha)\in W(m,m,0;\psi)$ must satisfy $\det X = 0$. In other words, $W(m,m,0;\psi)$ is contained in the hypersurface defined by the equation $\det X = 0$.

This logic can be generalized to cover the more general case $m+u \leq n$. For this we first introduce more notation. For each $m\times n$ matrix $X\in\R^{mn}$ with column vectors $\xx^{(1)},\dots,\xx^{(n)}$, let $\wtilde{X}$ be the $m\times (n-u)$ matrix with column vectors $\xx^{(u+1)},\dots,\xx^{(n)}$, and let $\wtilde\bfalpha = (\alpha_{u+1},\ldots,\alpha_n)$. Then let $\Gamma\subset\R^{mn}\times \R^n$ be the set of $(X, \bfalpha)\in\R^{mn}\times \R^n$ such that the determinant of every $m\times m$ minor of $\wtilde X\oplus\wtilde\bfalpha$ is $0$. Here $\wtilde X\oplus\wtilde\bfalpha$ is the $(m+1)\times(n-u)$ matrix whose first $m$ rows are the rows of $\wtilde X$ and whose last row is $\wtilde\bfalpha$.
%The dimension of this smaller set is $(m+1)n - 2(n-m-u+1)$

%
%\comdavid{Actually this definition of $\Gamma$ is incorrect. $\Gamma$ needs to be defined as the set of all $(X,\bfalpha)\in \mb R^{mn}\times \mb R^n$ such that the determinant of every $m\times m$ minor of $(\wtilde X,\wtilde\bfalpha)$ is $0$. The dimension of this smaller set is $(m+1)n - 2(n-m-u+1)$.}
%
%\commumtaz{I have removed the reference (Cf. \cite[Lemma 1]{DickinsonHussain}) from the lemma below as our proof is correcting a mistake in \cite{DickinsonHussain}.  I have put the reference in the remark below and removed it fromThe proof is}

\begin{lemma}
\label{dimlemma}
For $m+u\le n$, the set $W(m,n,u;\psi)$ is contained in $\Gamma$. Moreover, 
\[\dim \Gamma = \gamma \df (m+1)n - 2(n-m-u+1).\]

%\commumtaz{A bit confused that why $\dim\Gamma\leq \dim W$,  as it follows from the corollary below? Instead it should be other way around}
\end{lemma}
%[need to update proof to accommodate new statement\ddr]
\begin{proof}
For the first statement, fix $(X,\bfalpha)\in W(m,n,u;\psi)$ and let $Y$ be an $m\times m$ minor of $\wtilde X$. Then if $\bfbeta$ is the vector whose coordinates are the coordinates of $\bfalpha$ corresponding to the indices of the columns of $Y$, then it is easy to see that $(Y,\bfbeta) \in W(m,m,0;\psi)$. We previously showed that this implies $\det Y = 0$. Since $Y$ was arbitrary, we have $(X,\bfalpha) \in \Gamma$.

Since, $\Gamma\subset\R^{mn}\times \R^n$ is defined to be the set of $(X, \bfalpha)\in\R^{mn}\times \R^n$ such that the determinant of every $m\times m$ minor of $(\wtilde X,\wtilde\bfalpha)$ is $0$. The second statement readily follows on subtracting the number of dependent coordinates coming from the condition that the determinant of every $m\times m$ minor of $(\wtilde X,\wtilde\bfalpha)$ is $0$, from the full dimension i.e. $(m+1)n$.
\end{proof}

\begin{remark}
This Lemma is an analogue of \cite[Lemma 1]{DickinsonHussain} but we should point out that the calculation in \cite{DickinsonHussain} has a small error: if we consider the set of matrices $\wtilde X$ such that the $m-1$ certain prescribed columns of $\wtilde X$ are equal to certain prescribed values, then it is not true that the dimension of this set is always equal to the number claimed. It may be higher, e.g. if we prescribe that all the columns must equal zero. However, the claim is true in the generic case and this is why the computation still gives the correct answer. %Proving it is a little annoying though. 

\end{remark}

\begin{theorem}\label{thm3}
Let $m+u\leq n$ and let $\psi$ be an approximating function. Let $f$
be a dimension function such that the maps $r\mapsto r^{-m(n-m-u+1)}f(r)$ and $r^{-mn-m-u+1}f(r)$ are both monotonic. Let $U$ be a nonempty open set. Then
\[
\HH^f\big(W(m,n,u;\psi)\cap U\big)=
\begin{cases}
0 & {\rm {if}} \ \
\sum\limits_{q=1}^{\infty}f(\Psi
(q))\Psi(q)^{-(m-1)n-m-u+1} q^{m+u-1}<\infty\\
\HH^f(\Gamma\cap U) & {\rm {if}} \ \
\sum\limits_{q=1}^{\infty}f(\Psi
(q))\Psi(q)^{-(m-1)n-m-u+1} q^{m+u-1}=\infty.
\end{cases}
\]
\end{theorem}

Analogues of Corollaries \ref{cor2} and \ref{co1} may be stated in a similar way. For the sake of brevity we only state the analogue of Corollary \ref{cor3}.

\begin{corollary}
For $m+u\le n$ and $\tau > 0$. Then for every nonempty open set $U \subset \R^{mn}\times \R^n$ we have
\[
\dim\big(W(m,n,u;\tau)\cap U\big)= \left\{ \begin{array}{ll}
(m+1)n-2\left(n-\frac{m+u}{\tau+1}\right) &{\rm \ if}\quad \tau \geq \frac{1}{m+u-1},\\
&\\
(m+1)n - 2(n-m-u+1) & {\rm \ if}\quad \tau\le \frac{1}{m+u-1}.
\end{array}\right.
\]
\end{corollary}

\medskip

%END NEW%

Similar to the classical setting, next we consider the variant of the set ${\rm Bad}(m,n)$ complementary to $W(m,n,u;\psi)$. To be precise, for each $0\leq u \leq n$ let $\Bu$ denote the set of all $( X, \bfalpha) \in
\R^{mn}\times \R^n$ for which there exists a constant $C(X,\bfalpha)$ such that
\begin{multline*}\label{mixedbad}
\max\{|\qq\cdot
\xx^{(1)}-\alpha_1-p_1|,\ldots,|\qq\cdot \xx^{(u)}-\alpha_u-p_u|, \\ |\qq\cdot
\xx^{(u+1)}-\alpha_{u+1}|, \ldots,|\qq\cdot
\xx^{(n)}-\alpha_n|\}\geq C(X,\bfalpha) \cdot |\qq|^{-\frac{m+u}{n}+1}
\end{multline*}
%\comdavid{changed exponent $-(m+1)/n$ to $-\frac{m+u}{n}+1$}
\noindent for all integer vectors
$(\pp, \qq)\in\Z^u\times\Z^m$.
It is an easy consequence of Theorem \ref{thm2} that for $m+u>n$, $\Bu$ is a null-set, i.e. $|\Bu|_{(m+1)n}=0$. On the other hand, by Lemma \ref{dimlemma}, for $m+u\leq n$ we have $\R^{(m+1)n} \setminus \Gamma \subset \Bu$ and thus $|\Bu|_{(m+1)n}=1$. This raises the natural question of the Hausdorff dimension of $\Bu$ whenever $m+u>n$, and of $\Bu\cap\Gamma$ whenever $m+u \leq n$, as well as of the associated fibers $${\rm Bad}_\bfalpha(m,n,u) = \{X\in \R^{mn} : (X,\bfalpha)\in \Bu\}.$$ The set  ${\rm Bad}_\0(m,1,0)$ was first studied in \cite{Hussain}, wherein it was proved to have maximal dimension. Later in \cite{HussainKristensen}, this result was extended to general linear forms ${\rm Bad}_\0(m,n,u)$ allowing $u \neq 0$. %Note that by \cite[Corollary 7.12]{Falconer_book}, this implies that $\Bu$ has full dimension. 
We strengthen this result as follows:

%NEW%
%In this paper, we prove that the set $\Bu$ is hyperplane winning, a notion which we will be defined in section \ref{sec:BadProof}. For the time being, we note that the property of being hyperplane winning implies the property of being winning for Schmidt's game which in turn implies the property of having full Hausdorff dimension. Thus, our result is the following:

\begin{theorem}\label{bad}
The set $\Bu$ is hyperplane winning, and in particular has full Hausdorff dimension. If $m +u\leq n$, then the set $\Bu\cap \Gamma$ is hyperplane winning relative to $\Gamma$.
\end{theorem}

It is worth noting that the method of proof of the present paper is immediately applicable to the situation when $u=0$, and in that sense we are significantly strengthening the result of \cite{HussainKristensen}. As stated before, the hyperplane winning property passes automatically to games played on certain fractals, therefore, this theorem implies that, for any fractal set $K$, the set $\Bu\cap K$ is hyperplane winning. This further implies that $\Bu$ has full Hausdorff dimension when intersected with the fractal set $K$.

\medskip

\section{Preliminaries and auxiliary results}

To kick off this section we first define the basic concepts of Hausdorff measure and dimension.

\subsection{Hausdorff Measure and Dimension}
\label{hm}

Below is a brief introduction to Hausdorff $f$--measure and dimension. For further details see \cite{BernikDodson, Falconer_book}. Let $F\subset \R^n$. For any $\rho>0$, a countable collection $\{B_i\}$ of balls in $\R^n$ with diameters $\mathrm{diam} (B_i)\le \rho$ such that $F\subset \bigcup_i B_i$ is called a $\rho$--cover of $F$. Define for a right continuous, monotonically increasing function $f: \R_{\ge 0} \rightarrow \R_{\ge 0}$ with $f(t) > 0$ for $t > 0$,
\[
\HH_\rho^f(F)=\inf \sum_if(\mathrm{diam}(B_i)),
\]
where the infimum is taken over all possible $\rho$--covers of $F$. The function $f$ is called a dimension function, and the Hausdorff $f$--measure of $F$ is
\[
\HH^f(F)=\lim_{\rho\to 0}\HH_\rho^f(F).
\]
In the particular case where $f(r)=r^s$ with $ s>0$, we write $\HH^s$ for $\HH^f$ and the measure is referred to as $s$--dimensional Hausdorff measure. The Hausdorff dimension of $F$ is denoted by $\dim F $ and is defined as
\[
\dim F :=\inf\{s\in \R^+\;:\; \HH^s(F)=0\}.
\]

\subsection{Slicing}\label{sslicing}

We now state a result due to Beresnevich and Velani \cite{BeresnevichVelani4}, referred to as the \emph{slicing lemma}, which is a key ingredient in the proof of Theorem \ref{thm2}. We include the result mainly for completeness, as its application is identical to the one in \cite{HussainKristensen2}. Before we state the result it is necessary to introduce a little notation.

Suppose that $V$ is a linear subspace of $\R^k$, $V^{\perp}$ will be used to denote the linear subspace of $\R^k$ orthogonal to $V$. Further $V+a:=\left\{v+a:v\in V\right\}$ for $a\in V^{\perp}$.

\medskip

\begin{lemma}[Slicing Lemma, {\cite[Lemma 4]{BeresnevichVelani4}}]\label{slicing} Let $l, k \in \N$ be such that $l\leq k$ and let $f$ be a dimension function such that $g(r) \df r^{-l} f(r)$ is also a dimension function. Let $B\subset \R^k$ be a Borel set and let $V$ be a $(k-l)$--dimensional linear subspace of $\R^k$. If for a subset $S$ of $V^{\perp}$ of positive $\HH^l$ measure$$\HH^{g}\left(B\cap(V+b)\right)=\infty \quad\forall \; b\in S,$$ \noindent then $\HH^{f}(B)=\infty.$
\end{lemma}

\subsection{Mass Transference Principle}\label{mtp}
We describe the mass transference principle for linear forms tailored for our use. The actual framework is broad ranging and deals with the limsup sets defined by a sequence of neighborhoods of `approximating' planes. The mass transference principle for linear forms naturally enables us to generalize the Lebesgue measure statements of linear forms to the Hausdorff measure statements. In its original form, it was derived by Beresnevich and Velani from the mass transference principle for simultaneous approximation \cite{BeresnevichVelani} using their `slicing' technique introduced in \cite {BeresnevichVelani4} and described above. %\comdavid{This used to say ``the one dimensional mass transference principle'' but there is nothing one-dimensional about simultaneous approximation, which is covered in their original paper.}

Let $\xx^{(j)}$ denote the $j$th column vector of $X$. For $\qq\in \Z^m\setminus \{\0\}$, and $\pp\in\Z^n$,  the resonant set $R_{\pp, \qq}$ is defined by
%\begin{equation*}\label{resonantset}
%R_{\pp, \qq}:= \bigcup_{\pp} \left\{ (X, \bfalpha)\in
%\R^{mn}\times \R^n:\qq X=\pp+\bfalpha \right\} =\bigcup_{\pp}R_{{\qq, \pp}}
%\end{equation*}
%\comdavid{This seems wrong, what if $\pp\in\Z^n\setminus(\qq\R^{mn})$ is close enough to $(\qq\R^{mn})$ that its neighborhood intersects it?}
%\commumtaz{here I was meaning to consider all those $\pp$ which are $|\pp|\leq |\qq|$}

%\noindent where for each $\pp$ %\in (\qq\R^{mn})\cap\Z^n$,
\[R_{ \pp, \qq}=\left\{ (X, \bfalpha)\in
\R^{mn}\times \R^n:\qq X-\bfalpha =\pp\right\}= R_{p_1, \qq}\times\cdots\times R_{p_n, \qq}
\]
\noindent and \[R_{p_j,\qq}=\left\{ (\xx^{(j)},  \alpha_j):\qq\cdot\xx^{(j)}-\alpha_j =p_j\right\}.\]

It is then clear that the resonant sets are the affine subspaces of dimension $mn$, codimension $n$ and are contained in $W(m,n;\psi)$ for all functions $\psi$.

Let $\RR =\left\{ R_{\pp, \qq }:\qq \in \Z^m\setminus\{\0\}, \pp\in\Z^n\right\}$. Given an approximating function $\psi$ and a resonant set $R_{\pp, \qq }$, define the $\Psi$--neighbourhood of $R_{\pp,\qq}$ as
\begin{equation*}
\Delta\left( R_{\pp,\qq},\Psi (\vert\qq%
\vert )\right) =\left\{ (X, \bfalpha)\in
\R^{mn}\times \R^n: \text{dist}\left(
(X, \bfalpha),R_{\pp,\qq}\right) \leq \frac{\psi (\left\vert \qq\right\vert )}{%
\left\vert \qq\right\vert }\right\},
\end{equation*}
where $\mathrm{dist}(A, B):=\inf \{|a-b|:a\in A, b\in B\}.$ Notice that if $m=1$ then the resonant sets are points and the sets $\Delta\left( R_{\pp, \qq},\Psi (|\qq|)\right)$ are balls centred at these points. %This is of little concern to us, as we assume throughout that $m>1$, but it is nevertheless instructive.

Let
\[
\Lambda(m,n;\Psi)
=\{(X, \bfalpha)\in \R^{mn}\times \R^n:(X, \bfalpha)\in \Delta\left( R_{\pp, \qq},\Psi (|\qq%
| )\right) \text{for i.m. } (\pp, \qq)\in \Z^n\times\Z^m\setminus\{\0\}\}
\]
and define
\begin{equation*}
\label{limsup1}
\Delta_{\qq}(\Psi):= \bigcup_{|\pp|\leq |\qq|}
\Delta( R_{\pp, \qq}
,\Psi(|\qq|)).
\end{equation*}
%where
%\[J=\left\{(\pp, \qq)\in\Z^n\times\Z^m\setminus\{\mathbf{0}\}:|\pp|\leq |\qq|\right\}.\]
Then, $\Lambda(m,n; \Psi)$ can be written as a limsup set so that
\begin{equation*}
\label{limsup2} \Lambda(m,n;\Psi) =
\bigcap_{N=1}^{\infty}\bigcup_{|\qq|=N}^\infty{}\Delta_{\qq}( \Psi).
\end{equation*}

\medskip

\begin{theorem}[Mass Transference Principle, \cite{AllenBeresnevich}]\label{MTPlinearforms}
Let $\RR$ and $\Psi$ as above be given. Let $f$ be a dimension function such that $g(r):=r^{-mn}f(r)$ is also a dimension function and such that the map $r\mapsto r^{-(m+1)n}f(r)$ is monotonic. Suppose for any ball $B$ in $\R^{(m+1)n}$
\begin{equation*}
\HH^{(m+1)n}\left(B\cap\Lambda(m, n;g(\Psi)^{\frac{1}{n}})\right)=
\HH^{(m+1)n}(B).
\end{equation*}
Then for any ball $B\in \R^{(m+1)n}$
\begin{equation*}
\HH^f\left(B\cap\Lambda(m,n;\Psi)\right)=
\HH^f(B).
\end{equation*}
\end{theorem}

This theorem was originally proved with an additional hypothesis on the system of subspaces $\RR$, see \cite[Theorem 3]{BeresnevichVelani4}. However, we need the stronger version in \cite{AllenBeresnevich} in order to prove Theorems \ref{thm2} and \ref{thm1}.

\subsection{An interlude on classical Diophantine approximation}

Below is the statement of a generalized form of a Khintchine--Groshev-type theorem for $W(m,n;\psi)$ which to our knowledge is a new result and interesting in its own right. The reader will notice that it is a special case of Theorem \ref{thm2}, namely the case $u=n$, which correspond to the nearest integers lying in a full lattice and so the classical Diophantine approximation. We include it here, as it is a key ingredient in the proof of Theorem \ref{thm2}. 
%
%\comdavid{Previously there was ``and as it can be easily proved by appealing to the mass transference principle stated above.'', but I found this very confusing (both because it's not clear why more justification is needed for including the theorem, and because the proof of the theorem is given on the next page and involves more than just directly applying the mass distribution principle) so I removed it.}

\begin{theorem}\label{thm1}

Let $\psi$ be an approximating function. Let $f$ and $r\mapsto r^{-mn}f(r)$ be dimension functions such that $r\mapsto r^{-(m+1)n}f(r)$ is monotonic. Then
\[\HH^{f}(W(m,n;\psi))
= \left\{
\begin{array}{cl}
0 &
{\rm \ if} \ \sum
\limits_{q=1}^{\infty}f(\Psi (q))\Psi(q)^{-mn}q^{m+n-1}<\infty,\\
&\\
\HH^f(\R^{mn}\times \R^n)& {\rm \ if} \ \sum
\limits_{q=1}^{\infty}f(\Psi (q))\Psi(q)^{-mn}q^{m+n-1}=\infty.
\end{array}\right.\]
\end{theorem}

Note that even in the one dimensional setting there is no condition of monotonicity imposed on the approximating function, in contrast to the case of estimating the dimension of a fiber corresponding to a fixed value of $\alpha$ (either zero or nonzero), where the monotonicity condition on the approximating function cannot be removed due to the Duffin--Schaeffer counterexample, see \cite {BDV, DuffinSchaeffer}.

Note that if the dimension function $f$ is such that $r^{-(m+1)n}f(r)\to \infty$ as $r\to 0$ then $\HH^f(\R^{mn}\times \R^n)=\infty$ and Theorem \ref{thm1} is the analogue of the classical result of Jarn\'{\i}k (see \cite{DickinsonVelani, Jarnik3}). In the case where $f(r) := r^{(m+1)n}$ the Hausdorff measure $\HH^f$ is proportional to the standard $(m+1)n$--dimensional Lebesgue measure supported on $\R^{(m+1)n}$ and the result is the natural analogue of the Khintchine--Groshev theorem for $ W(m,n;\psi)$. A singly metric analogue of Theorem \ref{thm1} can be found in Bugeaud's paper \cite{Bugeaud5}.

We state the following special case of Theorem \ref{thm1}, which is a consequence of \cite[Theorem 15]{Sprindzuk}.

%\comdavid{Removed the word ``corollary'' as to me this suggests a result that can be proven by first proving the main theorem and then deducing the corollary, not vice-versa.} \comdavid{Also, in the third paragraph of the introduction it is stated that the doubly metric Khinchin--Groshev can be found in Schmidt or Sprind\v zuk, not Cassels. Though I wasn't able to find it precisely in any of those places...}

\begin{theorem}[{\cite[Theorem 15]{Sprindzuk}}]\label{cor1} \noindent Let $\psi$ be an approximating function. Then

\[| W(m,n;\psi) |_{(m+1)n} = \left\{
\begin{array}{cl}
0 &
{\rm \ if} \
\sum\limits_{q=1}^\infty{}\psi(q)^n q^{m-1}<\infty,\\
&\\
1& {\rm \ if} \
\sum\limits_{q=1}^\infty{}\psi(q)^n q^{m-1}=\infty.
\end{array}\right.\]
\end{theorem}

%In fact, we will deduce the divergence case of Theorem \ref{thm1} from the divergence  case of Theorem \ref{cor1}. It is very similar to the result of Schmidt \cite{Schmidt5}, but in a sense more general as there is no assumption of monotonicity on the approximating function. 
For the doubly metric case, the Hausdorff dimension for the set $W(m,n;\psi)$ was established by Dodson in \cite{Dodson} and, for the singly metric case, by Levesley in \cite {Levesley}. A slightly more general form of Dodson's result was established by Dickinson in \cite{Dickinson3}.

\subsection{Proof of Theorem \ref{thm1}}

%\begin{proof}[Proof of Theorem \ref{thm1}]{\mbox{ }}

\subsubsection{The convergence case} Notice that the set $W(m,n;\psi)$ can be written in the following limsup form:
\begin{equation*}
W(m,n;\psi) =\underset{N=1}{\overset{\infty }{%
{\bigcap }}}\,{\bigcup_{h=N}^\infty }\,\underset{\left\vert \qq%
\right\vert =h}{{\bigcup }}\Delta\left( R_{\pp, \qq},\Psi (\vert\qq%
\vert )\right).
\end{equation*}
For each resonant set $R_{\pp, \qq}$
the set $\Delta(R_{\pp, \qq},\Psi(|\qq|))$
can be covered by a collection of $(mn+n)$-dimensional closed hypercubes with disjoint interior and sidelength comparable with $\Psi(|\qq|)$. It can be readily verified that the number $C$ of such hypercubes satisfies
\[
C\lessless \Psi(|\qq|)^{-mn}|\qq|^n.
\]
Thus, $W(m,n;\psi)$ can be written as the limsup of a sequence of hypercubes whose total ``$f$-dimensional cost'' is
\begin{eqnarray*}
\sum_{h = 1}^\infty \sum_{|\qq|=h}
C f(\Psi(|\qq| ))
&\lessless&
\sum_{h = 1}^\infty \sum_{|\qq|=h}
\Psi(|\qq|)^{-mn}|\qq|^n f(\Psi(|\qq| )) \\
&\lessless &\sum_{h = 1}^\infty h^{m+n-1}f(\Psi (h))\Psi (h)
^{-mn} < \infty.
\end{eqnarray*}
Thus by the Hausdorff--Cantelli lemma \cite[Lemma 3.10]{BernikDodson},
we have $\HH^f(W(m,n;\psi))=0$, as required. 
%\comdavid{The previous version of this proof relied on the incorrect inequality $\mcal H^f(\textup{set}) \leq \textup{cost}(\textup{cover})$. This inequality is incorrect because the cost of a $\rho$-cover only gives an upper bound for $\mcal H^f_\rho$, not for $\mcal H^f$. I find that the easiest way to do things properly is via the Hausdorff--Cantelli lemma, so I've rewritten the proof to use that.}

\subsubsection{The divergence case}

The divergence case is an easy consequence of the mass transference principle discussed in the section \6 \ref{mtp}. In view of this we shall use the divergence part of Theorem \ref{cor1} (Cassels' theorem) and the mass transference principle to prove the divergence part of Theorem \ref{thm1}. With reference to the framework of \6\ref{mtp} we note that $$\Lambda(m, n; \frac{1}{m+1}\Psi) \subset W(m,n;\psi) \subset \Lambda(m, n; \Psi)$$ 
%the set $\Lambda(m,n;\Psi)$ coincides with the set $W(m,n;\psi)$. \comdavid{This is not true, you have to do something like $\Lambda(c\Psi) \subset W(m,n;\psi) \subset \Lambda(C\Psi)$.} 
Hence the divergence case follows.
%\end{proof}

\section{Proof of Theorem \ref{thm2}}

\subsection{The convergence case} The convergence part,  although similar to the convergence part of Theorem \ref{thm1}, is included for completeness. 

For $\qq\in \Z^m\setminus \{\0\}$, and $\pp_u=(p_1,\ldots, p_u)\in\Z^u$,  the resonant set $R_{\pp_u, \qq}$ is defined by

\[R_{ \pp_u, \qq}=\left\{ (X, \bfalpha)\in
\R^{mn}\times \R^n:\qq X-\bfalpha =\pp_u\right\}= R_{p_1, \qq}\times\cdots\times R_{p_u, \qq}\times \cdots\times R_{p_n, \qq},
\]
\noindent where \[R_{p_j,\qq}=\left\{ (\xx^{(j)},  \alpha_j):\qq\cdot\xx^{(j)}-\alpha_j =p_j\right\} \quad (1\leq j\leq u)\]
\noindent and \[R_{p_i,\qq}=\left\{ (\xx^{(i)},  \alpha_j):\qq\cdot\xx^{(i)}-\alpha_i =0\right\} \quad (u+1\leq i\leq n). \]

It is then clear that the resonant sets are the affine subspaces of dimension $mn$, codimension $n$ and are contained in $W(m,n, u;\psi)$ for all functions $\psi$.
Further, the set $W(m,n, u;\psi)$ can be written in the following limsup form:
\begin{equation*}
W(m,n, u;\psi) =\underset{N=1}{\overset{\infty }{%
{\bigcap }}}\,{\bigcup_{h=N}^\infty }\,\underset{\left\vert \qq%
\right\vert =h}{{\bigcup }}\Delta\left( R_{\pp_u, \qq},\Psi (\vert\qq%
\vert )\right),
\end{equation*}
where

\begin{equation*}
\Delta\left( R_{\pp_u,\qq},\Psi (\vert\qq%
\vert )\right) =\left\{ (X, \bfalpha)\in
\R^{mn}\times \R^n: \text{dist}\left(
(X, \bfalpha),R_{\pp_u,\qq}\right) \leq \frac{\psi (\left\vert \qq\right\vert )}{%
\left\vert \qq\right\vert }\right\},
\end{equation*}
is the   $\Psi$--neighbourhood of $R_{\pp_u,\qq}$ for any resonant set $R_{\pp_u, \qq }$.
For each resonant set $R_{\pp_u, \qq}$
the set $\Delta(R_{\pp_u, \qq},\Psi(|\qq|))$
can be covered by a collection of $(mn+n)$-dimensional closed hypercubes with disjoint interior and sidelength comparable with $\Psi(|\qq|)$. It can be readily verified that the number $C$ of such hypercubes satisfies
\[
C\lessless \Psi(|\qq|)^{-mn}|\qq|^u.
\]
Thus, $W(m,n, u;\psi)$ can be written as the limsup of a sequence of hypercubes whose total ``$f$-dimensional cost'' is
\begin{eqnarray*}
\sum_{h = 1}^\infty \sum_{|\qq|=h}
C f(\Psi(|\qq| ))
&\lessless&
\sum_{h = 1}^\infty \sum_{|\qq|=h}
\Psi(|\qq|)^{-mn}|\qq|^u f(\Psi(|\qq| )) \\
&\lessless &\sum_{h = 1}^\infty h^{m+u-1}f(\Psi (h))\Psi (h)
^{-mn} < \infty.
\end{eqnarray*}
Thus by the Hausdorff--Cantelli lemma, we have $\HH^f(W(m,n, u;\psi))=0$, as required. 

\subsection{The divergence case}

As discussed earlier the statement of the Theorem essentially reduces to two
cases, the finite measure case where $r^{-(m+1)n}f(r)\to C>0$ as
$r\to 0$ and the infinite measure case in which $r^{-(m+1)n}f(r)\to \infty$ as $r\to 0$. Therefore, we split the proof of
the Theorem \ref{thm2} into two parts, the finite measure case and
the infinite measure case.

Before proceeding, we will need the following key lemma, which will
make our proofs work.

\begin{lemma}

\label{lem:bi-lip}

Let $S \subset \mat_{(m+u-n)\times n}(\mathbb{R}) \times \mathbb{R}^n$ of full Lebesgue measure. Let $A \subset \GL_{n\times n}(\mathbb{R})$ be a subset of a subspace of dimension $(n-u)n$ such that $A$ has positive $(n-u)n$--dimensional Lebesgue measure. Then, the set
\begin{equation*}
\Lambda = \left\{\left(
\begin{pmatrix}
X \\
YX
\end{pmatrix},\bfalpha X\right) \in \mat_{(m+u) \times n}(\mathbb{R}) \times \mathbb{R}^n: X \in A, (Y, \bfalpha) \in S \right\}
\end{equation*}
has full Lebesgue measure inside $A \times S$, with the product structure being the one implicit in the affine system.
\end{lemma}
\begin{proof}

By Fubini's theorem, it suffices to show that for almost all $X\in A$, the fiber $$\Lambda_X = \{(YX,\bfalpha X) : (Y,\bfalpha) \in S\}$$ has full measure. In fact, this is true for all $X\in \GL_{n\times n}(\R)$; this follows from the invariance of Lebesgue measure under right multiplication by the invertible matrix $X$.
\end{proof}
\begin{remark}
For a detailed proof for the analogous homogeneous situation we refer to Lemma $4.1$ of \cite{HussainKristensen2}.  
\end{remark}

\subsection{Finite measure case}
As in \cite{HussainKristensen2}, in order to proceed, we will make some restrictions. Let $\epsilon > 0$ and $N > 0$ be fixed but arbitrary. Define
\begin{equation*}
A_{\epsilon, N} = \left\{(X, \bfalpha) \in \mat_{m\times n}(\mathbb{R})\times \R^n : \epsilon <
\det(\wtilde{X}) < \epsilon^{-1}, \quad \max_{1 \leq i,j\leq n}
\abs{x_{ij}} \leq N \right\}
\end{equation*}
where $\wtilde{X}$ denotes the $n \times n$-matrix formed by the first $n$ rows of the matrix $X$. The set is of positive measure for $\epsilon$ small enough and $N$ large enough, and as $\epsilon$ decreases and $N$ increases, the set fills up $\mat_{(m+1)\times n}(\mathbb{R})$ with the exception of the null-set of matrices $X$ such that $\wtilde{X}$ is singular.

We will prove that the divergence assumption implies that the set $W(m,n,u;\psi)$ is full in $A_{\epsilon, N}$. We will translate the statements about $W(m,n,u;\psi)$ into ones about usual Diophantine approximation, specifically regarding sets of the form $W_\bfalpha\left( m+u-n,n;c\psi \right)$, $c > 0$.

Consider the set of $n$ affine forms in $m+u-n$ variables defined by the pair $(\hat{X}, \hat\bfalpha) \in \R^{(m+u-n)n}\times\R^n$. Suppose furthermore that these linear forms satisfy the inequalities
\begin{equation}
\label{eq:5}
\norm{\mathbf{r} \hat{X}-\hat{\bfalpha}}_i \leq
\frac{\psi(\abs{\mathbf{r}})}{nN}, \quad 1 \leq i \leq n,
\end{equation}
for infinitely many $\rr \in \mathbb{Z}^{m+u-n} \setminus \{\0\}$, where $\norm{\mathbf{x}}_i$ denotes the distance from the $i$th coordinate of $\mathbf{x}$ to the nearest integer. A special case of Theorem \ref{cor1} states that the divergence condition of our theorem implies that the set of such pairs $(\hat{X}, \hat\bfalpha)$ is full in $\mat_{(m+u-n)\times n}(\R)\times \R^n$, and hence in particular also in the image of $A_{\epsilon, N}$ under the map $$\left(\begin{pmatrix}
\wtilde{X} \\
\hat{X}\wtilde{X}
\end{pmatrix}, \hat{\bfalpha}\wtilde{X}\right)\mapsto (\hat{X}, \hat{\bfalpha}).$$

Now, suppose that $(X, \bfalpha) \in A_{\epsilon, N}$ is such that $(\hat{X}, \hat\bfalpha)$ is in the set defined by \eqref{eq:5}, where $(\hat{X}, \hat\bfalpha)$ is the uniquely determined pair such that 
$$
(X, \bfalpha) = \left(\begin{pmatrix}
\wtilde{X} \\
\hat{X}\wtilde{X}
\end{pmatrix}, \hat{\bfalpha}\wtilde{X}\right).
$$
We claim that $(X,\bfalpha)$ is in $W(m,n,u;\psi)$. Indeed, let $(\mathbf{r}_k)_k$ be an infinite sequence of integer vector such that the inequalities \eqref{eq:5} are satisfied for each $k$, and let $\mathbf{p}_k$ be the nearest integer vector to $\mathbf{r}_k \hat{X} - \bfalpha$. Now define $\mathbf{q}_k = (\mathbf{p}_k, \mathbf{r}_k)$. The inequalities defining $W(m,n,u;\psi)$ will be satisfied for these values of $\mathbf{q}_k$, since
\begin{eqnarray*}
\abs{\mathbf{q}_k \begin{pmatrix}
I_u & 0 \\
X_u & X^\prime
\end{pmatrix}-\bfalpha} &\asymp& \abs{\mathbf{q}_k
\begin{pmatrix}
I_n \\
\hat{X}
\end{pmatrix} \wtilde{X}-\bfalpha}\\ & \asymp_\epsilon & \abs{\left(\pm \norm{\mathbf{r}_k\cdot \xx^{(1)}+\hat{\bfalpha}_1}_1, \dots,
\pm \norm{\mathbf{r}_k\cdot \xx^{(n)}+\hat{\bfalpha}_n}_n\right) \wtilde{X}},
\end{eqnarray*}
where $\hat{X}$ splits into $n$ column vectors $\left(\xx^{(1)}, \ldots, \xx^{(n)}\right)$, $X=(X_u \ X^\prime)$.  Where the implied  constant may depend upon $m, n, \epsilon$ but not on $\qq_k$. The $i$th coordinate of the first vector is at most $\psi(\abs{ q})/nN$, so carrying out the matrix multiplication, using the triangle inequality and the fact that $\abs{x_{ij}} \leq N$ for $1 \leq i, j \leq n$ shows that
\begin{equation*}
\abs{\mathbf{q}_k \begin{pmatrix}
I_u & 0 \\
X_u & X^\prime
\end{pmatrix}-\bfalpha}_i < \psi(\abs{\mathbf{q}_k}).
\end{equation*}
Applying Lemma \ref{lem:bi-lip}, the divergence part of Theorem
\ref{thm2} follows in the case of Lebesgue measure.

\subsection{Infinite measure case}
The infinite measure case of the Theorem \ref{thm2} can be easily
deduced from the following lemma.

\begin{lemma}
\label{lem:infmeslem}
Let $\psi$ be an approximating function and let $f$ be a dimension function such that $g(r) \df r^{-n(n-u)}f(r)$ and $r \mapsto r^{-(m+u-n)n}g(r)$ are also dimension functions, $r\mapsto r^{-(m+u+1-n)n}g(r)$ is monotonic, and $r^{-(m+1)n}f(r)\to\infty$ as $r\to{}0$. If
\[
\sum \limits_{q=1}^{\infty}f(\Psi(q))\Psi(q)^{-mn}q^{m+u-1}=\infty,
\]
then
\[
\HH^{f}(W(m,n,u;\psi))=\infty.
\]
\end{lemma}

%\comdavid{I decided to skip to section 5 at this point, I will read through this proof and the next section on a later iteration.}
%
%\commumtaz{Please read through the proof, I have made some changes to the previous version, mostly notational. }
\begin{proof}
The proof relies on a bi-Lipschitz injective map which preserves measure, just like the proof of Lemma $4.2$ of \cite{HussainKristensen2}. As in the finite measure case, we fix $\epsilon > 0$, $N \geq 1$ and let
\begin{equation*}
A_{\epsilon, N} = \left\{(X, \bfalpha) \in \mat_{m\times n}(\mathbb{R})\times \R^n : \epsilon <
\det(\wtilde{X}) < \epsilon^{-1}, \quad \max_{1 \leq i,j\leq n}
\abs{x_{ij}} \leq N \right\},
\end{equation*}
where $\wtilde{X}$ denotes the $n \times n$-matrix formed by the first $n$ rows of $X$. We also define the set
\begin{equation*}
\wtilde{A}_{\epsilon, N} = \left\{\wtilde{X} \in \GL_{n}(\mathbb{R}) :
\epsilon < \det(\wtilde{X}) < \epsilon^{-1}, \quad \max_{1 \leq i,j\leq n} \abs{x_{ij}} \leq N \right\}.
\end{equation*}
Consider a subspace
\begin{equation*}
A_0 =\left\{X\in \wtilde{A}_{\epsilon, N}: X=\begin{pmatrix}
I_u & 0 \\
X_u & X^\prime
\end{pmatrix}\right\},
\end{equation*}
of $\wtilde{A}_{\epsilon, N}$, which is of positive $(n-u)n-$dimensional Hausdorff measure. For an appropriately chosen constant $c>0$ depending only on $m,n,\epsilon$ and $N$, we find that the map
\begin{eqnarray*}
% \label{eq:2}
\eta : W(m+u-n, n, c \psi) \times A_0 &\rightarrow&
W(m,n,u;\psi), \\ ((Y, \bfalpha), X)  &\mapsto&
\left(\begin{pmatrix}
X \\
YX
\end{pmatrix}, \bfalpha X\right)
\end{eqnarray*}
is a Lipschitz embedding. Indeed, it is evidently injective as $X$ is invertible for all the domain. We have
\begin{eqnarray*} \HH^{f}(W(m,n,u;\psi))
&\geq&\HH^f\left(\eta\left(W(m+u-n, n, c \psi)\times A_0\right)\right) \notag \\
&\asymp&\HH^{f}\left(W(m+u-n, n, c \psi)\times A_0\right).\label{sec:infinite-measure}
\end{eqnarray*}
Now after taking into account the conditions on the dimension functions as stated in the lemma along with the divergence part of Theorem \ref{thm1} for $W(m+u-n, n; c \psi)$ and $f=g$, we apply the slicing lemma \ref{slicing} across $A_0$ to get the infinite measure under the appropriate assumption of divergence. %The details are identical to those in \cite{HussainKristensen2}, and we omit them.
\end{proof}

\section{Proof of Theorem \ref{thm3}}\label{dsmuleqn}
As in the preceding results, the convergence case is a rather dull affair, and boils down to constructing an appropriate fine cover and applying the Hausdorff--Cantelli lemma \cite[Lemma 3.10]{BernikDodson}. We omit the details and proceed with the divergence case. % \comdavid{Changed ``Borel--Cantelli'' to ``Hausdorff--Cantelli'' and reworded appropriately}

Let $p = m+u-1$. For clarity we only discuss the case $m+u=n$, or equivalently $n-p=1$. Define the map
\[
\eta:\R^{mp}\times \R^p\times \R^{m-1} \to \Gamma
\]
via the formula
\begin{multline*}
\eta\left((\xx^{(1)},\dots,\xx^{(p)}),(\alpha_1,\ldots,\alpha_p)
,(b_{1},\dots,b_{m-1})\right)\\ =
\left(\Big(\xx^{(1)},\dots,\xx^{(p)}, \sum_{j=1}^{m-1}b_j\xx^{(u+j)}\Big),\Big(\alpha_1,\ldots,\alpha_p,\sum_{j=1}^{m-1} b_j \alpha_{u+j}\Big)\right).
\end{multline*}
For the more general case $m+u \leq n$, we need a slightly more complicated expression for $\eta$ with more linearly dependent columns in the matrix. However, the ideas should be clear from the above simpler case.

Fix $C \geq 1$,
\begin{equation}
\label{inclusion}
\left((\xx^{(1)},\dots, \xx^{(p)}),(\alpha_1,\ldots,\alpha_p) \right)\in
W(m,m+u-1,u;C^{-1}\psi),
\end{equation}
and $(b_1,\ldots,b_{m-1})\in\R^{m-1}$ such that $\sum_{j=1}^{m-1} |b_j| \leq C$. Note that by the triangle inequality,
\begin{eqnarray*}
\left|\qq\cdot\sum_{j=1}^{m-1}b_{j}\xx^{(u+j)}-\sum_{j=1}^{m-1} b_j \alpha_{u+j}\right|
&\leq & \sum_{j=1}^{m-1} |b_j| \cdot \big| \qq \cdot \xx^{(u+j)} - \alpha_{u+j}\big|\\
&\leq & \left(\sum_{j=1}^{m-1}|b_{j}|\right) C^{-1} \psi(|\qq|)\leq
\psi(|\qq|),
\end{eqnarray*}
for the infinitely many integer vectors $\qq$ we obtain from condition \eqref{inclusion}. It follows that
\[
\eta\Big(W(m,m+u-1,u;C^{-1}\psi)\times \left[-\tfrac{C}{m-1},\tfrac{C}{m-1}\right]^{m-1}\Big) \subset W(m,n,u;\psi).
\]
In the more general case $m+u\leq n$, we have
\begin{equation}
\label{etageneral}
\eta\Big(W(m,m+u-1,u;C^{-1}\psi)\times \left[-\tfrac{C}{m-1},\tfrac{C}{m-1}\right]^{(m-1)(n-p)}\Big) \subset W(m,n,u;\psi).
\end{equation}
To proceed further, we will use the fact that $\eta$ is a locally bi-Lipschitz embedding when restricted to the set $\Gamma'\times \R^{(m-1)(n-p)}$, where
\[
\Gamma' = \{(X,\bfalpha) : \rk(\wtilde X\oplus\wtilde\bfalpha) = m-1\}.
\]
Note that $\Gamma'$ is relatively open and dense in $\Gamma$.

We will need the following lemma \cite[Proposition 2.2]{HSS}.
\begin{lemma}\label{fulllemma}
Let $f$ be a dimension function, let $L\subset\R^l$, and let $\eta:L\to \R^k$ be a bi-Lipschitz embedding. Then $\HH^f(L)\asymp\HH^f(\eta(L))$.
\end{lemma}

We are now in a position to prove the infinite measure case of Theorem \ref{thm3}.

\begin{lemma}
\label{lemmixedinf}
Let $m+u \leq n$. Let $\psi$ be an approximating function and let $f$ be a dimension function such that $g(r) = r^{-(m-1)(n-p)}f(r)$ is also a dimension function, and such that
\[
\lim_{r\to 0} r^{-\gamma} f(r) = \infty.
\]
Let $U \subset \R^{mn}\times \R^n$ be a non empty open set. If
\[
\sum_{q=1}^{\infty} f(\Psi(q)) \Psi(q)^{-(m-1)n-p} q^p=\infty,
\]
then
\[
\HH^{f}(W(m,n,u;\psi)\cap U)=\infty = \HH^f(U).
\]
\end{lemma}
\begin{proof}
Let $V \subset \R^{mp}\times \R^p$ be a nonempty open set. Since
\[
\lim_{r\to 0} r^{-(m+1)p} g(r) = \infty,
\]
we have $\HH^g(V) = \infty$. On the other hand, since
\[
\sum_{q=1}^\infty g(\Psi(q)) \Psi(q)^{-mp} q^p=\infty,
\]
by Theorem \ref{thm2}, for every $C\geq 1$ we have
\[
\HH^g\big(W(m,m+u-1,u;C^{-1}\psi)\cap V\big) = \HH^g(V) = \infty.
\]
Thus by \cite[Corollary 7.12]{Falconer_book}, for every nonempty open set $V' \subset \R^{m-1}$ we have
\[
\HH^f\big((W(m,m+u-1,u;C^{-1} \psi)\cap V) \times V'\big) = \infty.
\]
Now choose $V$, $V'$, and $C$ such that
\[
\eta(V\times V') \subset U, \;\;\;\;\; V' \subset \left[\tfrac C{m-1},\tfrac C{m-1}\right]^{(m-1)p}
\]
and such that $\eta$ is bi-Lipschitz on $V\times V'$. Then by \eqref{etageneral} and Lemma \ref{fulllemma},
\[
\HH^f\big(W(m,n,u;\psi)\big) = \infty.
\]

Now suppose that $\liminf_{r\to 0} r^{-\gamma} f(r) < \infty$. Then either $\HH^f(\Gamma) = 0$, or $\HH^f$ is proportional to $\HH^\gamma$ on $\Gamma$. Thus for the remainder of the proof, we assume that $f(r) = r^\gamma$.

As before, let $g(r) = r^{-(m-1)(n-p)} f(r)$, and suppose that
\[
\sum_{q=1}^{\infty} f(\Psi(q)) \Psi(q)^{-(m-1)n-p} q^p = \sum_{q=1}^\infty g(\Psi(q)) \Psi(q)^{-mp} q^p = \infty.
\]
Again, by Theorem \ref{thm2}, for every $C\geq 1$ and for every open set $V \subset \R^{mp}\times \R^p$ we have
\[
\HH^g\big(W(m,m+u-1,u;C^{-1}\psi)\cap V\big) = \HH^g(V)
\]
Since $g(r) = r^{(m+1)p}$, this means that $W(m,m+u-1,u;C^{-1}\psi)$ is of full Lebesgue measure in $\R^{mp}\times \R^p$. Thus,
\[
W(m,m+u-1,u;C^{-1}\psi) \times \R^{(m-1)(n-p)}
\]
is of full measure in $\R^{mp}\times \R^p\times \R^{(m-1)(n-p)}$. Since $\eta$ is smooth, by \eqref{etageneral} and Sard's theorem \cite [Theorem II.3.1]{Sternberg},  the set $W(m,n,u;\psi)$ is of full measure in the image of $\eta$, which is in turn of full measure in $\Gamma$. So $W(m,n,u;\psi)$ is of full measure in $\Gamma$.
\end{proof}

\section{Proof of Theorems  \ref{thm:BadHAW}} 

We begin with a description of the hyperplane game. 
\subsection{Hyperplane  Winning.} 
Fix $0 < \beta < 1$. The $\beta$-hyperplane game on $\R^k$ has two players, Alice and Bob, and is played as follows:
\begin{enumerate}
\item Bob chooses an initial closed ball $B_0 = B(\mathbf{x}_0, \rho_0)$.
\item After Bob's $n$th move $B_n$, Alice `deletes a neigbourhood of a hyperplane $A_n$', i.e. she picks an affine hyperplane $A_n \subset \R^k$, which must be avoided by Bob in his next move.
\item Bob now chooses a ball $B_{n+1} = B(\mathbf{x}_{n+1}, \rho_{n+1})$ satisfying
\begin{equation*}
B_{n+1} \subset B_n \setminus A_n^{(\beta\rho_n)} \text{ and } \rho_{n+1} \ge \beta \rho_n.
\end{equation*}
Here, $A_n^{(\beta\rho_n)}$ denotes a closed tubular neighbourhood of $A_n$ of radius $\beta\rho_n$. If this is not possible, Alice wins the game.
\end{enumerate}

A set $S$ is said to be \emph{$\beta$-hyperplane winning} if Alice has a strategy ensuring that \begin{equation*}
\bigcap_{n=1}^\infty B_n \cap S \neq \emptyset
\end{equation*}
If $S$ is $\beta$-hyperplane winning for all $\beta > 0$, it is said to be \emph{hyperplane winning}. 

The game is a modification of the classical Schmidt game \cite{Schmidt1}, which we will also need. in this game we have two parameters $\alpha> 0$ and $\beta>0$. In this game, instead of removing a neighbourhood of a hyperplane, Alice chooses a closed ball $W_n \subseteq B_n$ of radius $\alpha \rho_n$. Also, at each step Bob must choose a closed ball $B_{n+1} \subseteq W_n$ of radius $\rho_{n+1} = \beta \alpha \rho_n$. The conditions for winning the game are unchanged. A set $S$ is $(\alpha, \beta)$-winning if Alice has a winning strategy for this set of parameters. It is $\alpha$-winning if it is $(\alpha, \beta)$-winning for any $\beta > 0$.

%It is this latter property which we will prove for the set $\Bu$.

%As in the preceding cases, we will transfer results from classical Diophantine approximation to the absolute value setting. In this case, the result needed appears not to be published anywhere, so we will prove it here. It is a strengthening of the result of Kleinbock \cite{Klenibock4}, which makes explicit our remark in the introduction that Einsiedler and Tseng \cite{EinTse} considerably strengthen this result.
%
%\begin{theorem}
%\label{thm:BadHAW}
%The set $\B$ is hyperplane winning.
%
%\end{theorem}

\subsection{Proof of Theorem \ref{thm:BadHAW}}

The idea is to modify the proof that $\Bbo$ is hyperplane winning found in \cite{BFS1}. Using the ideas of \cite{EinsiedlerTseng}, this approach can be extended to prove that $\Balpha$ is hyperplane winning. This was stated as a result in a preprint of \cite{BFS1}, though it did not make it to the final paper. In fact, the approach is sufficiently flexible to allow us to modify the proof to show that $\Badmn$ is hyperplane winning.

Before starting the technical proof, let us explain the strategy in a little more detail. For the classical Schmidt game, Einsiedler-Tseng \cite{EinsiedlerTseng} modified Schmidt's \cite{Schmidt2} original proof that $\Bbo$ is winning to show that $\Balpha$ is winning for any $\alpha$. Proceeding in the same way, but with the hyperplane game in place of the classical game, one can extend the proof in \cite{BFS1} to show that $\Balpha$ is hyperplane winning. The strategy is devised so that from the $k$th step of the game, Bob is only able to chose from a set of affine forms which satisfy the inequality defining $\Balpha$ for all $\qq$ with $\vert \qq \vert$ bounded by an explicit parameter of the proof. We will provide a proof of this, based on the strategy devised in \cite{BFS1} with modifications according to the ideas of \cite{EinsiedlerTseng}. 

The key new insight is that in her $k$th move, Alice does not need to know the exact value of $\bfalpha$, but only to know the value to within an accuracy of
\[
c \rho_k^{\frac{m}{m+n}},
\]
which allows us to complete the proof of Theorem \ref{thm:BadHAW}.

We now proceed with an explicit description of Alice's strategy. We will only describe the steps in the construction of the strategy and refer the reader to \cite{BFS1} for the full details. Initially, we describe how she plays in the case of $\Bbo$. This is the strategy shown to be hyperplane winning in \cite{BFS1}. Given a matrix $X \in \R^{mn}$, we append a standard unit vector to each column, i.e. we consider the following matrix,
\[
\widetilde{X} = \begin{pmatrix}
X \\ I_n
\end{pmatrix}.
\]
Note that $X \in \R^{mn}$ is an element of $\Bbo$ if ad only if there is a constant $c > 0$ such that for any $(\qq,\pp) \in \Z^m \times \Z^n$,
$$
\vert \xx \vert^n \cdot  \vert(\qq, \pp) \widetilde{X}\vert^n > c.
$$ 

Let $\lambda = n/(m+n)$, $R>1$ and $\delta = R^{n(n+m)^2}$. The strategy devised in \cite{BFS1} ensures that that for any $i \in \N$, the system of inequalities
\begin{equation}
\label{eq:game1}
\begin{split}
0 < \vert \qq \vert &< \delta R^{m(\lambda+i)} \\
\vert (\qq,\pp)\widetilde{X} \vert &< \delta R^{-n(\lambda + i) - m},
\end{split}
\end{equation}
has no solutions $(\qq, \pp) \in \Z^{m+n} \setminus \{\0\}$. It is easily seen that this gives a winning strategy for $\Bbo$.

Unfortunately, the strategy is not easily described, as one needs to consider simultaneously a dual setting. In this setup, one considers instead the matrices
\[
\widetilde{X^T} = \begin{pmatrix}
X^T \\ I_m
\end{pmatrix},
\]
and introduce a new paramter, $\delta^T = R^{-m(m+n)^2}$. The condition which is ensured simultaneously to \eqref{eq:game1} is that for any $j \in \N$, 
the system of inequalities
\begin{equation}
\label{eq:game2}
\begin{split}
0 < \vert \pp \vert &< \delta^T R^{n(\lambda+j)} \\
\vert (\pp,\qq)\widetilde{X^T} \vert &< \delta R^{-m(\lambda + i) - n},
\end{split}
\end{equation}
has no solutions $(\pp, \qq) \in \Z^{n+m} \setminus \{\0\}$.

At the heart of constructing a winning strategy is the following lemma \cite[Lemma 5.3]{BFS1}.

\begin{lemma}
\label{lem:Bad0game}
For $R \in \R$ sufficiently large, the following holds
\begin{enumerate}[i)]
\item When the game is at stage $k_i$, suppose that for any $X \in B(X_{k_i}, \rho_{k_i})$, the system of inequalities \eqref{eq:game1} has no  solutions $(\qq, \pp) \in \Z^{m+n} \setminus \{\0\}$; and the system of inequalities \eqref{eq:game2}  has no solutions $(\pp, \qq) \in \Z^{n+m} \setminus \{\0\}$ with $j=i-1$. Then, Alice has a strategy such that after finitely many stages at stage $h_i$, say, the system of inequalities \eqref{eq:game2}  has no solutions $(\pp, \qq) \in \Z^{n+m} \setminus \{\0\}$ with $j=i$ for any $X \in B(X_{h_i}, \rho_{h_i})$.
\item Dually, at stage $h_j$, suppose that for any $X \in B(X_{h_i}, \rho_{h_i})$, the system of inequalities \eqref{eq:game1} has no  solutions $(\qq, \pp) \in \Z^{m+n} \setminus \{\0\}$; and the system of inequalities \eqref{eq:game2}  has no solutions $(\pp, \qq) \in \Z^{n+m} \setminus \{\0\}$ with $j=i$.  Then, Alice has a strategy such that after finitely many stages at stage $k_{j+1}$, say, the system of inequalities \eqref{eq:game1} has no solutions $(\qq, \pp) \in \Z^{m+n} \setminus \{\0\}$ with $i = j+1$ for any $X \in B(X_{k_{j+1}}, \rho_{k_{j+1}})$
\end{enumerate}
\end{lemma}

It clearly suffices to construct such a strategy, as one can then recursively construct a hyperplane winning strategy for $\Bbo$. This also explains why it is only necessary to know $\bfalpha$ to a certain precision in our later modification of the strategy, as we are successively avoiding more and more affine forms. We will be more precise later.

The proofs of the two parts of the lemma are each others dual, as should be apparent from the definitions given. We will comment on the first part, so suppose that for all $X \in B =B(X_{k_i}, \rho_{k_i})$, the non-existence of interger solutions to \eqref{eq:game2} with $j=i-1$ and \eqref{eq:game1} is ensured as in the statement of the lemma. Define the set
\[ 
S = \left\{(\pp,\qq) \in \Z^{n+m} : \text{ for some $X \in B$, \eqref{eq:game2} is satisfied with $j=i$}\right\}.
\]
The $\R$-span of $S$ has dimension at most $n$. If the dimension is less than $n$, we extend this subspace in some arbitrary way to a space of dimension equal to $n$. Now, pick an orthonormal basis $\mathcal{Y} = \{\yy_1, \dots, \yy_n\}$ for this subspace of $\R^{n+m}$. Let $\bb_u$ denote the $u$th column of $\widetilde{X^T}$ and define the matrix 
\[ 
M(X,\mathcal{Y}) = (\bb_u \cdot \yy_v)_{1 \le u,v \le n},
\]
and let $D$ denote the determinant of this matrix. It is now shown using Cramer's rule that under the given assumptions,
\begin{equation}
\label{eq:game3}
\vert D \vert \le n \sqrt{n} R^{-(m+n)(i+1)} \max\{\vert D_{11} \vert, \vert D_{12}\vert, \dots, \vert D_{nn} \vert\},
\end{equation}
The idea is now to construct a strategy for Alice such that after a finite number of steps, \eqref{eq:game3} cannot hold. In this way, she will have avoided all solutions to \eqref{eq:game2} with $j=i$, and the game can continue with the dual construction.

To prove that it is possible for Alice to ensure that \eqref{eq:game3} cannot hold after a finite number of steps, we will need additional notation. Let $v \in \{0,1,\dots, n\}$ and define the set
\[
\Omega_v = \{(I,J) \in \{0,1,\dots, n\} \times \{0,1,\dots, n\} : \#I = \#J = v\}.
\]
For each $\omega = (I,J) \in \Omega_v$, let $D_\omega(X) = \det(\bb_{i_k} \cdot \yy_{j_l})_{1 \le k,l \le v}$, where $I = \{i_1, \dots i_v\}$ and $J = \{j=1, \dots, j_v\}$, i.e. the minor of the determinant $D$ corresponding to the choices $I$ and $J$. Now, define the vector
\[
\overline{M}_v(X) = (D_\omega(X))_{\omega \in \Omega_v} \in \R^{\binom{n}{v}^2}.
\]
Finally, let $M_v(X) = \Vert \overline{M}_v(X) \Vert$, where $\Vert \cdot \Vert$ denotes the euclidean norm. For a ball $B \in \R^{mn}$, we let
\[
M_v(B) = \sup_{X \in B} M_v(X).
\]
Two cases require special treatment. For $v=0$, $\Omega_0 = \{(\emptyset, \emptyset)\}$, and we let $D_{(\emptyset, \emptyset)}(X) = 1$ for all $X$, so that $M_{0}(B) = 1$ identically. For $v = -1$, $\Omega_{-1} = \emptyset$, and $M_{-1}(B) = 0$ identically.

The final ingredient in the proof is Lemma 6.1 from \cite{BFS1}, which is now stated.

\begin{lemma}
\label{lem:finite_game}
Let $0 < \beta < 1/3$, $\sigma \in \R$ and $v \in \{0,1,\dots, n\}$. There is a constant 
\[ 
\nu_v = \nu_v(m,n,\beta, \sigma) > 0,
\]
such that for any $\mu_v \in (0,\nu_v]$ and any orthonormal system $\mathcal{Y} = \{\yy_1, \dots, \yy_n\} \subset \R^{m+n}$, Alice can win the finite game with rules
\begin{enumerate}[i)]
\item Bob chooses a closed ball $B \subseteq \R^{mn}$ of radius $\rho_B < 1$ and with $\max_{X \in B} \Vert X \Vert \le \sigma$.
\item Bob and Alice play according to the rules of the hyperplane game until the radius of Bob's ball $B_v$ is less than $\mu_v \rho_B$.
\item Alice wins if for any $X \in B_v$,
\[
M_v(X) > \nu_v \rho_B M_{v-1}(B_v).
\]
\end{enumerate}
\end{lemma}

The reader will note that the inequality of the lemma in the case $v= N$ essentially is of the form \eqref{eq:game3}, and indeed it is exactly what is required to obstruct existence of solutions. On choosing the parameter $R$ large enough, the previous lemma ensures the validity of Lemma \ref{lem:Bad0game}. %See \cite{BFS1} for the details.

The construction of the strategy is now based on some results of Schmidt \cite{Schmidt2}, which are also proved in \cite{BFS1}. We condense them into a single lemma.

\begin{lemma}
\label{lem:schmidt}
Fix $X \in B$ and $v \in \{0, \dots, n\}$. Then,
\begin{enumerate}[i)]
\item For any $\omega \in \Omega_v$, $\Vert \nabla D_\omega(X)\Vert_{\text{op}} \le v M_{v-1}(X)$.
\item For any $\omega \in \Omega_v$, $\Vert \nabla \nabla D_\omega(X)\Vert_{\text{op}} \le v^2 M_{v-2}(X)$.
\item There exist constants $\epsilon_1, \epsilon_2 > 0$ depending only on $m,n$ and $\sigma$, such that if $M_v(X) \le \epsilon_1 M_{v-1}(X)$, then 
\[
\max_{\omega \in \Omega_v} \Vert \nabla D_\omega(X) \Vert_{\text{op}}  > \epsilon_2 M_{v-1}(X).
\]
\end{enumerate}
\end{lemma}
Where $\Vert \cdot  \Vert_{\text{op}}$ stands for the operator norm.

We now proceed by induction in $v$. For $v=0$, the claim  is trivial, so suppose Lemma \ref{lem:finite_game} is true for $v-1$. In other words, for some $\nu_{v-1} > 0$ and some $\mu_{v-1} \in (0,\nu_{v-1}]$, if $B_{v-1}$ is the first ball chosen by Bob with radius less than $\mu_{v-1}\rho_B$, then for all $X \in B_{v-1}$,
\[
M_{v-1}(X) > \nu_{v-1} M_{v-2}(B_{v-1}).
\]

Applying inequality i) of lemma \ref{lem:schmidt} and using the mean value theorem, it is straightforward to show that 
\begin{equation}
\label{eq:game4}
M_{v-1}(B_{v-1}) \le v M_{v-1}(X).
\end{equation}

At this point, we are finally able to describe the strategy, Alice should follow to win the finite game. There are two cases. First, if $M_v(X) > \epsilon_1 M_{v-1}(X)$ for all $X \in B$, Alice will make arbitrary moves until $\rho(B_v) < \mu_v \rho_B$. This suffices, since \eqref{eq:game4} implies that
\[
M_v(X) \ge \frac{\epsilon_1}{v} M_{v-1}(B_{v-1}) \ge \frac{\epsilon_1}{v} M_{v-1}(B_{v}),
\]
since $B_v \subseteq B_{v-1}$. Hence, chosing $\nu_v < \frac{\epsilon_1}{v}$ suffices in this case.

Hence, we need only worry about the case $M_v(X) > \epsilon_1 M_{v-1}(X)$ for some $X \in B$. In this case, lemma \ref{lem:schmidt} iii) holds, so that
\[
\Vert \nabla D_\omega(X) \Vert_{\text{op}}  > \epsilon_2 M_{v-1}(X)
\]
for some $\omega \in \Omega_v$. We consider the linearization $L$ of  $D_\omega$ around $X$, 
\[
L(A') = D_\omega(X) + \nabla_{A' - A} D_\omega(A).
\]
 Let $\mathcal{L} = L^{-1}(0)$, an affine subspace of $\R^{mn}$. Alice removes the neighbourhood $\mathcal{L}^{(\beta \rho(B_{v-1}))}$ from $B_{v-1}$ and continues to make arbitrary moves until $\rho(B_v) < \mu_v \rho_B$. Because of lemma \ref{lem:schmidt} iii), if $A' \in B_{v-1} \setminus \mathcal{L}^{(\beta \rho(B_{v-1}))}$, $\vert L(A') \vert$ is fairly large. On the other hand,  lemma \ref{lem:schmidt} iii) allows us to control how much $L(A')$ deviates from $D_\omega(A')$. Putting all estimates together, one finds that chosing $\nu_v = \frac{1}{2v}\beta^2 \epsilon_2 \mu_{v-1}$ will ensure the conclusion. See \cite{BFS1} for the calculations. This concludes our treatment of $\Bbo$.

Now, to get from $\Bbo$ to $\Balpha$, we show how to modify the strategy using the ideas of \cite{EinsiedlerTseng}. In fact, their proof shows directly that $\Balpha$ is hyperplane winning. The added complication comes from the fact that they play their game on certain fractals. We describe the strategy in the following.

We will use their notation, which again involve certain matrices. For $X$ an $m\times n$-matrix and $\bfalpha \in \R^m$, define 
\[
L_X(\bfalpha) = \begin{pmatrix}
I_m & X & -\alpha \\
0 & I_n & 0 \\
0 & 0 & 1
\end{pmatrix}.
\]
The affine lattice associated to $X$ and $\bfalpha$ is given by
\[
L_X(\bfalpha) (\Z^{m+n} \times \{1\}) = \left\{ L_X(\bfalpha)\begin{pmatrix} \pp \\ \qq \\ 1\end{pmatrix} : \pp \in \Z^m, \qq \in \Z^n\right\} \subset \R^{m+n} \times \{1\}.
\]
Note that the last coordinate is always equal to $1$, so by abuse of notation, we will also write $L_X(\bfalpha) (\Z^{m+n})$ for the projection of the above set onto the first $m+n$ coordinates, and we will refer to this set as the affine lattice associated to $X$ and $\bfalpha$. 

The space of all affine unimodular lattices $\Omega_{m+n, {\rm aff}}$  in $\R^{m+n}$ consists of all translates $\Lambda + \bfbeta$, where $\Lambda = g \Z^{m+n}$ for some $g \in \SL(m+n, \R)$ is a unimodular lattice and $\bfbeta \in \R^{m+n}$. Note that the affine lattices associated to some pair $(X, \bfalpha)$ are of this form.

Finally, we will require a one-parameter action on the space $\Omega_{m+n, {\rm aff}}$. The appropriate one is the left action by matrices of the form
\[
g_t = \begin{pmatrix} e^{t/m} I_m &0 &0 \\ 0 & e^{-t/n}I_n & 0 \\ 0 & 0 & 1 \end{pmatrix},
\]
on $\R^{m+n+1}$ and in particular on the invariant subspace $\R^{m+n} \times \{1\}$ as above. Observe that the action defined on $\R^{m+n} \times \{1\}$ induces an action on $\Omega_{m+n, {\rm aff}}$.

The key idea is to apply the following result of Dani \cite[Theorem 2.20]{Dani}. For an affine lattice $\Gamma \in \Omega_{m+n, {\rm aff}}$, we let $\delta(\Gamma)$ denote the length of the shortest non-zero vector in $\Gamma$.
\begin{theorem}
\label{thm:Dani}
The pair $(X, \bfalpha)$ is in $\Balpha$ if and only if for some $c > 0$,
\[
\inf_{t \in \R} \delta(g_t L_X(\bfalpha) \Z^{m+n}) \ge c.
\]
\end{theorem}
 
 As we have seen, Alice has a $\beta$-winning strategy for the $\Bbo$-game. By \cite{BFKRW}, this implies that Alice has an $\alpha_0$-winning strategy for the usual Schmidt game for some $\alpha_0 > 0$. 
Suppose Bob has chosen the ball $B_\ell = B(A_\ell, \rho_\ell)$ at his last step of the hyperplane game. We will show how Alice should respond in her following two moves.
 
 Initially, let $\alpha = \delta \alpha_0$ and for $\beta > 0$ let $\beta_0 = \delta^{-1} \beta$, where $\delta> 0$ is chosen appropriately ($\delta = 2^n/\mu(B(0,1))$ will do, where $\mu$ is the $n$-dimensional Lebesgue measure, see \cite{EinsiedlerTseng}). Note that $\alpha_0 \beta_0 = \alpha \beta$. The usual Schmidt game is $(\alpha_0, \beta_0)$-winning. Hence, Alice may choose a ball $B_\ell = B(A_\ell', \alpha_0 \rho_\ell)$ according to this strategy. 
 
 Next, $t_\ell$ must be chosen so that for any $D \in \R^{mn}$, $g_{t_\ell} L_D(0) = L_{(\alpha_0 \rho_\ell)^{-1} D} (0) g_{t_\ell}$ for all $D$. This is done by setting
\[
\alpha_0 \rho_\ell = \exp\left(-\left(\tfrac1m + \tfrac1n\right) t_\ell\right).
\]
We then define the affine lattice $x_\ell$ by
\[
x_\ell = g_{t_\ell} L_{A_\ell'}(\bfalpha) \Z^k = L_\0(e^{\frac1n t_\ell}\bfalpha) g_{t_\ell} L_{A_\ell'}(\0) \Z^k.
\]

Now, let $\vv$ be a vector in $x_\ell$ of shortest norm. For a vector $\xx \in \R^{m+n} =  \R^m \times \R^n$, let $\xx_p$ denote projection onto the first $m$ components (the particle space in the language of \cite{EinsiedlerTseng}) and let $\xx_t$ denote projection onto the last $n$ coordinates  (the time space). There is a proper, affine subspace $\mathcal{L}$ of $\R^{mn}$ such that $(L_A(0) \vv)_p = 0$ if and only if $A \in \mathcal{L}$. In fact, just as above in the $\bfalpha = 0$ case, if $A \notin \mathcal{L}^{(\epsilon)}$, then 
\begin{equation}
\label{eq:game}
\Vert (L_A(0) \vv)_p \Vert \ge \epsilon \Vert v_t \Vert.
\end{equation}
Alice will remove this hyperplane for an appropriately chosen $\epsilon > 0$. The value $\epsilon = 2 \delta \alpha_0 \rho_\ell$ works. 

Note that in \cite{EinsiedlerTseng}, the authors proceed by letting Alice choose yet another ball $W_\ell$ inside $B'_\ell$ but disjoint from $\mathcal{L}^{(\epsilon)}$. This is due to the fact that they are concerned with the usual Schmidt game, but in the remainder of their argument, all that is really needed is \eqref{eq:game}. Hence, there is no problem in just removing the hyperplane neighbourhood and letting Bob pick the next ball.

In her next move, Alice plays according to a winning strategy for the hyperplane game with target set $\Bbo$ with the same parameter, i.e. $\beta = 2 \delta \alpha_0$. This is to ensure that however the game ends, the resulting point will be an element in $\Bbo$.

Let $X$ be the final point when the above strategy is followed. Fix $\gamma_1>0$ small enough that $g_t L_X(0) \Z^k$ does not contain a non-zero element of norm $\le \gamma_1$ by Theorem \ref{thm:Dani} with $\bfalpha = 0$. Then, fix $\gamma_2 > 0$ so that $L_A \Z^k$ does not contain any element with norm $\le \gamma_2$, which is possible since $\bfalpha \notin \Z^k$. Finally, pick $c > 0$ so that $\Vert D \vv_t \Vert \le c \Vert \vv_y \Vert$ for all $D \in B(0,1)$. In \cite{EinsiedlerTseng}, it is shown by induction that for any $\ell \ge 0$,
\[
B(0,r) \cap g_{t_\ell} L_X(\bfalpha) \Z^k = \emptyset
\]
for 
\[
r = \tfrac{1}{2} \min\{ \gamma_1, \gamma_2\} (1+c)^{-2} (\alpha\beta)^{n/(m+n)}.
\]
We will not repeat the argument here, but appealing to Theorem \ref{thm:Dani}, this immediately implies that $X \in \Balpha$, so that this set is hyperplane winning.

Finally, to complete the proof, note that in the above, when Alice is responding to Bob's ball $B_k = B(X_k,\rho_k)$, she only needs know the value of $\bfalpha$ to within an accuracy of
\[
c \rho_k^{\frac{m}{m+n}},
\]
where $c > 0$ is a small constant. Using the claim, we complete the proof as follows. Suppose that in the hyperplane game with target set ${\rm Bad}(m,n)$, Bob plays the ball $B_k = B((X_k,\bfalpha_k),\rho_k)$. Then Alice can interpret this move as corresponding to the move $B_k' = B(X_k,\rho_k)$ in the hyperplane game with target set $\Balpha$, where $\bfalpha \in B(\bfalpha_k,\rho_k)$ is an unknown parameter. Since $\rho_k < c \rho_k^{\frac{m}{m+n}}$ for all sufficiently large $k$, we see that $\bfalpha$ is known to within an accuracy of $c \rho_k^{\frac{m}{m+n}}$ which is sufficient for Alice to use the strategy described above to generate a hyperplane-neighborhood $A_k' = \mathcal L^{(\beta\rho_k)}$ with which to respond to Bob's move $B_k'$. In the hyperplane game with target set ${\rm Bad}(m,n)$, Alice plays the corresponding hyperplane-neighborhood
\[
A_k = A_k' \times \R^n = \{(X,\bfalpha) : X\in A_k', \; \bfalpha\in \R^n\}.
\]
The game continues with an eventual result of $(X,\bfalpha)$ in the game with target set ${\rm Bad}(m,n)$ and an eventual result of $X$ in the game with target set $\Balpha$. Since $X\in \Balpha$ if and only if $(X,\bfalpha)\in {\rm Bad}(m,n)$, this shows that ${\rm Bad}(m,n)$ is hyperplane winning.

\qed

\section{Proof of Theorem \ref{bad}}\label{sec:BadProof}
We now proceed to deduce Theorem \ref{bad}  from Theorem \ref{thm:BadHAW} by transporting a winning strategy from the classical setup to the new one. This is accomplished by appealing to the following result. 

\begin{proposition}
\label{prop:bad_transfer}
Let $k, l \in \N$, let $\mathcal{U} \subset \R^k$ be a nonempty open set whose complement is an algebraic variety, and let $f : \mathcal{U} \rightarrow \R^l$ be a differentially surjective map, locally $C^1$ conjugate to a projection $\pi$ on $\mathcal{U}$. Then the preimage under $f$ of any hyperplane winning set in $\R^l$ is hyperplane winning.
\end{proposition}

The proof of the Proposition depends on two older results, which we re-state here for easy reference. The first is \cite[Lemma 3.9]{FKMS}.

\begin{lemma}
\label{lem:FKMS}
Let $k,D \in \N$ and $\beta \in (0,1]$.There exists a $\gamma > 0$, such that for any non-zero polynomial $f : \R^k \rightarrow \R$ of degree at most $D$, if Bob and Alice play the $\beta$-hyperplane game and if Bob's first move is $B_0 = B(0,1)$, then Alice has a strategy ensuring that the first ball chosen by Bob of radius less than $\gamma$ is disjoint from the set 
\[
\left\{\xx \in \R^k : \dist(\xx, f^{-1}(0)) < \gamma\right\}
\]
\end{lemma}

The other key ingredient is a special case of \cite[Theorem 2.4]{BFKRW}, which we state as a lemma.

\begin{lemma}
\label{lem:BFKRW}
Let $S \subseteq \R^k$ be hyperplane winning, $\UU \subseteq \R^k$ be open and let $f : \UU \rightarrow \R^k$ be $C^1$ non-singular map. Then, $f^{-1}(S) \cup \UU^c$ is hyperplane winning.
\end{lemma}

\begin{proof}
Alice begins by ensuring that we need not worry about the variety $\mathcal{U}^c$. This is accomplished by Lemma \ref{lem:FKMS}. Indeed, the variety $\mathcal{U}^c$ is the zero set of some polynomial map, so by the Lemma, Alice has a strategy ensuring that after finitely many steps, Bob must choose balls at a positive distance $\gamma > 0$ from this zero set. Hence, there is no loss of generality in assuming that the first ball chosen by Bob is contained in $\mathcal{U}$.

By Lemma \ref{lem:BFKRW}, applying a $C^1$ conjugation to a hyperplane winning set does not affect this property. Hence, in order to prove the proposition, it suffices to show that the preimage under a projection $\pi$ of a hyperplane winning set is itself hyperplane winning. This is however easy. Any ball $B(x,\rho)$ chosen by Bob in the preimage is sent to the ball $B(\pi(x), \rho)$ in the image. Here, Alice plays according to the winning strategy for this move and picks a hyperplane $A$ with neighbourhood $A^{(\beta \rho)}$, which is mapped to the hyperplane-neighbourhood $\pi^{-1}(A)^{(\beta \rho)}$ in the preimage. Evidently, this is a hyperplane winning strategy.
\end{proof}

With some additional care, one can extend Proposition \ref{prop:bad_transfer} to the following more general form, which may have other applications. We will not pursue this in the present paper, as we will only need the simpler version.

\begin{proposition}
\label{prop:bad_transfer_general}
Let $\mathcal{M}_1, \mathcal{M}_2$ be manifolds, let $\mathcal{U} \subset \mathcal{M}_1$ be a nonempty open set whose complement is an analytic variety and let $f : \mathcal{U} \rightarrow \mathcal{M}_2$ be a differentially surjective map, locally $C^1$ conjugate to a projection $\pi$ on $\mathcal{U}$. The preimage under $f$ of any hyperplane winning set in $\mathcal{M}_2$ is hyperplane winning.
\end{proposition}

In order to complete the proof of Theorem \ref{bad}, we need only ensure that the conditions of Proposition \ref{prop:bad_transfer} are satisfied. Before doing this we note that the defining inequalities take a particularly pleasing form. Namely, $(X,\bfalpha) \in \Bu$ if and only if
\begin{equation}
\label{eq:3}
\abs{(\pp,\qq)
\begin{pmatrix}
I_u & 0 \\
X_u & \wtilde{X}
\end{pmatrix}-\bfalpha}
\geq C(X,\bfalpha) \abs{\qq}^{-\frac{m+u}{n}+1}
\end{equation}
for all $(\pp,\qq) \in \Z^u \times \Z^m \setminus \{\0\}$, where $X$ is the matrix $(X_u, \ \wtilde{X})$.

We now split the proof into two parts depending upon the choices of $m$, $u$, and $n$.

\subsection{The case $m+u\leq n$} First, by way of motivating our proof when $m+u \leq n$, we discuss the case where $u=0$ and $m=n$. In this case, we consider the inequalities
\begin{equation*}
\abs{\mathbf{q}X-\bfalpha} \geq C(X)
\end{equation*}
for all but finitely many $\mathbf{q} \in \mathbb{Z}^m$. It is readily verified that unless there is linear dependence among the columns of $X$, this is trivially satisfied. Hence, in this simple subcase,
\[
\left(\R^{m^2}\setminus \{X\in\R^{m^2}:
\xx^{(1)},\ldots,\xx^{(m)} \text{ are linearly dependent}\} \right)\times \R^m \
\subset {\rm Bad}(m, m, 0).
\]
In other words, the set $ {\rm Bad}(m, m, 0)$ contains the set of matrices of full rank, which is hyperplane winning by \cite[Lemma 3.9]{FKMS}. This proves our main theorem in the particular case $m=n, u=0$.

Note that in fact we get the stronger inequality
\begin{equation}
\label{eq:2}
|\qq X-\bfalpha|>C(X)|\qq|\;\;\;\forall\qq\in \Z^{m}
\setminus \{\0\}
\end{equation}
in the set considered. This follows as an invertible matrix can only distort the unit ball by a specified amount. This is much stronger than the defining inequality of the set $ {\rm Bad}(m, m, 0)$ in this special case. This feature also applies to the more general setting when $m+u \leq n$ and further underlines the qualitative difference between the case $m+u \leq n$ and the converse $m+u > n$.

We now give a full proof in the case $m + u \leq n$. We will argue much in the spirit of the above. Fix $X \in \R^{mn}$ such that the matrix $\wtilde{X}$ in \eqref{eq:3} has full rank. By \cite[Lemma 3.9]{FKMS}, the set of such $X$ is hyperplane winning. Performing Gaussian elimination on the columns of a matrix of the form of \eqref{eq:3} implies the existence of an invertible $(n \times n)$-matrix $E(X)$ such that
\begin{equation*}
\begin{pmatrix}
I_u & 0 \\
X_u & \wtilde{X}
\end{pmatrix} =
\begin{pmatrix}
I_u & 0 & 0 \\
\hat{X} & I_m & 0
\end{pmatrix} E(X).
\end{equation*}
Applying this matrix from the right to a vector $(\pp, \qq)$, we see that
\begin{equation*}
(\pp,\qq)
\begin{pmatrix}
I_u & 0 \\
X_u & \wtilde{X}
\end{pmatrix} =
(\pp,\qq)
\begin{pmatrix}
I_u & 0 & 0 \\
\hat{X} & I_m & 0
\end{pmatrix} E(X) =
\begin{pmatrix}
\pp + \qq\hat{X} \\
\qq \\
\0
\end{pmatrix}^T E(X).
\end{equation*}

Multiplication by the matrix $E$ on the right hand side only serves to distort the unit ball in the absolute value to a different parallelepiped depending on $X$. This induces a different norm on the image, but by equivalence of norms on Euclidean spaces, this distortion can be absorbed in a positive constant. In other words,
\begin{equation}
\label{eq:1}
\abs{\left(\mathbf{p, q}\right)\left(
\begin{array}{cc}
I_u & 0 \\
X_u & \wtilde{X}
\end{array}
\right)-\bfalpha} \geq C(X)\abs{
\begin{pmatrix}
\pp + \qq\hat{X} \\
\qq \\
\0
\end{pmatrix}^T-\wtilde{\bfalpha}^T}
\end{equation}
for all $\left(\mathbf{p, q}\right)\in \Z^u \times \Z^{m} \setminus \{\0\}$, where $\wtilde{\bfalpha} = \bfalpha E(X)^{-1}$.

Finally, since the norm is the \emph{supremum} norm, we get
\begin{equation*}
\abs{
\begin{pmatrix}
\pp + \qq\hat{X} \\
\qq \\
\0
\end{pmatrix}^T-\wtilde{\bfalpha}^T} \geq C(\wtilde{\bfalpha})\abs{\qq},
\end{equation*}
for sufficiently large $\qq$. By \eqref{eq:1}, almost every $(X, \bfalpha)$ is in $\Bu$, and in fact with the stronger requirement from \eqref{eq:2}.  This completes the proof of Theorem \ref{bad} in the case where $m+u \leq n$.

\subsection{The case $m+u>n$}

In this case, the measure of the set $\Bu$ is zero, and the proof is not as elementary as above. Nevertheless, appealing to Proposition \ref{prop:bad_transfer} and taking a cue from the elementary proof, it is easily accomplished as follows.

%\comdavid{What used to be here didn't make much sense, so I tried to write something similar that made more sense.}

For each $(X,\bfalpha)\in \R^{mn}\times\R^n$, we can write
\[
X = \begin{pmatrix}
X_u & \wtilde X
\end{pmatrix}
= \begin{pmatrix}
X' & X'' \\
X''' & X''''
\end{pmatrix}
\]
where the columns are split into groups of size $u$ and $(n-u)$, and the rows are split into groups of size $(n-u)$ and $(m+u-n)$. Note that $X''$ is a square matrix. Now let $\mathcal U$ be the set of pairs $(X,\bfalpha)$ such that $\det X'' \neq 0$. Define the map $f:\mathcal U \to \mathcal M_2 \df \R^{(m+u-n)n}\times\R^n$ by letting
\[
f(X,\bfalpha) = \left(\begin{pmatrix}
X''' & X''''
\end{pmatrix}
\begin{pmatrix}
I_u & 0 \\
X' & X''
\end{pmatrix}^{-1},\bfalpha \begin{pmatrix}
I_u & 0 \\
X' & X''
\end{pmatrix}^{-1}\right).
\]
Now it is easy to see that $f$ is differentially surjective, and by Theorem \ref{thm:BadHAW}, the set ${\rm Bad}(m+u-n, n) \subset \mathcal{M}_2$ is hyperplane winning. Hence, to conclude by Proposition \ref{prop:bad_transfer}, it suffices to prove that the preimage of this set under $f$ is contained in $\Bu$.

However, this is easy. Fix $(X,\bfalpha)\in \R^{mn}\times\R^n$, and suppose that $(Y,\bfbeta) \df f(X,\bfalpha) \in {\rm Bad}(m+u-n, n)$. Then there is a constant $C(Y,\bfbeta) > 0$ such that for any $\mathbf{r} \in \mathbb{Z}^n$ and $\mathbf{s} \in \mathbb{Z}^{m+u-n}\setminus\{\0\}$,
\begin{equation}
\label{pr}
\abs{(\mathbf{r}, \mathbf{s})
\begin{pmatrix}
I_n \\
Y
\end{pmatrix} - \bfbeta} \geq C(Y,\bfbeta) \abs{\mathbf{s}}^{-\frac{m+u}{n} + 1}.
\end{equation}
On the other hand, if we let
\[
Z = \begin{pmatrix}
I_u & 0 \\
X' & X''
\end{pmatrix} \in \GL_n(\R)
\]
then
\[
\left(\begin{pmatrix}
I_u & 0 \\
X_u & \wtilde X
\end{pmatrix},\bfalpha\right)
= \left(\begin{pmatrix}
I_n\\
Y
\end{pmatrix}Z,\bfbeta Z\right)
\]
and thus for $\mathbf{p} \in \mathbb{Z}^u$ and $\mathbf{q} \in \mathbb{Z}^m \setminus \{\0\}$,
\begin{equation*}
\abs{(\mathbf{p}, \mathbf{q})
\begin{pmatrix}
I_u & 0 \\
X_u & \wtilde X
\end{pmatrix} - \bfalpha} = \abs{\left((\mathbf{p}, \mathbf{q})
\begin{pmatrix}
I_n \\
Y
\end{pmatrix} - \bfbeta\right) Z}
\geq C(X,\bfalpha) \abs{\mathbf{q}}^{-\frac{m+u}{n} +
1},
\end{equation*}
with the last inequality following by multiplying \eqref{pr} by the invertible matrix $Z$, and then using the fact that if $(\pp,\qq) = (\mbf r,\mbf s)$, then $|\qq| \geq |\mbf s|$. Thus \eqref{eq:3} is satisfied and $(X,\bfalpha) \in \Bu$.

%\commumtaz{I will update the references in the next iteration i.e. when the other parts of the paper are nearly complete}

\medskip

\noindent \emph{Acknowledgments.} The first-named author is supported by La Trobe University's start-up grant, the second-named author is supported by the Danish Research Council for Independent Research and the third-named author was supported by the EPSRC Programme Grant EP/J018260/1, and is now supported by a fellowship from the Royal Society. The authors would like to thank Stephen Harrap for good discussions and excellent company. Finally, the authors would like to thank the anonymous referees for their detailed comments which has significantly improved the quality of the paper. % Demi Allen

\providecommand{\bysame}{\leavevmode\hbox to3em{\hrulefill}\thinspace}
\providecommand{\MR}{\relax\ifhmode\unskip\space\fi MR }
\providecommand{\MRhref}[2]{%
  \href{http://www.ams.org/mathscinet-getitem?mr=#1}{#2}
}
\providecommand{\href}[2]{#2}

%
%\bibliographystyle{amsplain}
%
%\bibliography{bibliography}

\begin{thebibliography}{10}

\bibitem{AllenBeresnevich}
Demi Allen and Victor Beresnevich, \emph{A mass transference principle for
  systems of linear forms and its applications}, Compos. Math. \textbf{154}
  (2018), no.~5, 1014--1047. \MR{3798593}

\bibitem{Arnold}
V.~I. Arnol'd, \emph{Geometrical methods in the theory of ordinary differential
  equations}, Grundlehren der Mathematischen Wissenschaften [Fundamental
  Principles of Mathematical Science], vol. 250, Springer-Verlag, New
  York-Berlin, 1983, Translated from the Russian by Joseph Sz\"ucs, Translation
  edited by Mark Levi. \MR{695786}

\bibitem{BDV}
Victor Beresnevich, Detta Dickinson, and Sanju Velani, \emph{Measure theoretic
  laws for lim sup sets}, Mem. Amer. Math. Soc. \textbf{179} (2006), no. 846,
  x+91 pp. \MR{2184760}

\bibitem{BeresnevichVelani}
Victor Beresnevich and Sanju Velani, \emph{A mass transference principle and
  the {D}uffin-{S}chaeffer conjecture for {H}ausdorff measures}, Ann. of Math.
  (2) \textbf{164} (2006), no. 3, 971--992. \MR{2259250}

\bibitem{BeresnevichVelani4}
Victor Beresnevich and Sanju Velani, \emph{Schmidt's theorem, {H}ausdorff
  measures, and slicing}, Int. Math. Res. Not. (2006), Art. ID 48794, 24.
  \MR{2264714}

\bibitem{BernikDodson}
Vasilii Bernik and Maurice Dodson, \emph{Metric {D}iophantine approximation on
  manifolds}, Cambridge Tracts in Mathematics, vol. 137, Cambridge University
  Press, Cambridge, 1999.\MR{1727177}

\bibitem{BFKRW}
Ryan Broderick, Lior Fishman, Dmitry Kleinbock, Asaf Reich, and Barak Weiss,
  \emph{The set of badly approximable vectors is strongly {$C^1$}
  incompressible}, Math. Proc. Cambridge Philos. Soc. \textbf{153} (2012),
  no.~02, 319--339. \MR{2981929}

\bibitem{BFS1}
Ryan Broderick, Lior Fishman, and David Simmons, \emph{Badly approximable
  systems of affine forms and incompressibility on fractals}, J. Number Theory
  \textbf{133} (2013), no. 7, 2186--2205. \MR{3035957}

\bibitem{Bugeaud5}
Yann Bugeaud, \emph{An inhomogeneous {J}arn\'\i k theorem}, J. Anal. Math.
  \textbf{92} (2004), 327--349. \MR{2072751}

\bibitem{Cassels}
John W.~S. Cassels, \emph{An introduction to {D}iophantine approximation},
  Cambridge Tracts in Mathematics and Mathematical Physics, No. 45, Cambridge
  University Press, New York, 1957.\MR{0157947}

\bibitem{Dani}
S.~G. Dani,\emph{Divergent trajectories of flows on homogeneous spaces and Diophantine approximation}, J. Reine Angew. Math. \textbf{359} (1985), 55--89.\MR{0794799}

\bibitem{DaniLauNog}
S.~G. Dani, Michel Laurent, and Arnaldo Nogueira, \emph{Multi-dimensional
  metric approximation by primitive points}, Math. Z. \textbf{279} (2015),
  no.~3-4, 1081--1101. \MR{3318261}

\bibitem{DGY1}
D.~Dickinson, T.~Gramchev, and M.~Yoshino, \emph{First order pseudodifferential
  operators on the torus: normal forms, {D}iophantine phenomena and global
  hypoellipticity}, Proceedings of the {C}onference ``{D}ifferential
  {E}quations'' ({I}talian) ({F}errara, 1996), vol.~41, 1996, pp.~51--64
  (1997). \MR{1471014}

\bibitem{DGY2}
Detta Dickinson, Todor Gramchev, and Masafumi Yoshino, \emph{Perturbations of
  vector fields on tori: resonant normal forms and {D}iophantine phenomena},
  Proc. Edinb. Math. Soc. (2) \textbf{45} (2002), no.~3, 731--759. \MR{1933753}

\bibitem{DickinsonHussain}
Detta Dickinson and Mumtaz Hussain, \emph{The metric theory of mixed type
  linear forms}, Int. J. Number Theory \textbf{9} (2013), no.~1, 77--90.
  \MR{2997491}

\bibitem{DickinsonVelani}
Detta Dickinson and Sanju Velani, \emph{Hausdorff measure and linear forms}, J.
  Reine Angew. Math. \textbf{490} (1997), 1--36. \MR{1468922}

\bibitem{Dickinson}
H.~Dickinson, \emph{The {H}ausdorff dimension of systems of simultaneously
  small linear forms}, Mathematika \textbf{40} (1993), no.~2, 367--374.
  \MR{1260899}

\bibitem{Dickinson3}
\bysame, \emph{A note on the theorem of {J}arn\'\i k-{B}esicovitch}, Glasgow
  Math. J. \textbf{39} (1997), no.~2, 233--236. \MR{1460640}

\bibitem{Dodson}
M.~M. Dodson, \emph{A note on metric inhomogeneous {D}iophantine
  approximation}, J. Austral. Math. Soc. Ser. A \textbf{62} (1997), no.~2,
  175--185. \MR{1433207}

\bibitem{DRV_KAM}
M.~M. Dodson, J.~P\"oschel, B.~P. Rynne, and J.~A.~G. Vickers, \emph{The
  {H}ausdorff dimension of small divisors for lower-dimensional {KAM}-tori},
  Proc. Roy. Soc. London Ser. A \textbf{439} (1992), no.~1906, 359--371.
  \MR{1193007}

\bibitem{DRV4}
M.~M. Dodson, B.~P. Rynne, and J.~A.~G. Vickers, \emph{Diophantine
  approximation and a lower bound for {H}ausdorff dimension}, Mathematika
  \textbf{37} (1990), no.~1, 59--73. \MR{1067887}

\bibitem{DuffinSchaeffer}
R.~J. Duffin and A.~C. Schaeffer, \emph{Khintchine's problem in metric
  {D}iophantine approximation}, Duke Math. J. \textbf{8} (1941), 243--255.
  \MR{0004859}

\bibitem{EinsiedlerTseng}
Manfred Einsiedler and Jimmy Tseng, \emph{Badly approximable systems of affine
  forms, fractals, and {S}chmidt games}, J. Reine Angew. Math. \textbf{660}
  (2011), 83--97. \MR{2855820}

\bibitem{Falconer_book}
Kenneth Falconer, \emph{Fractal geometry: {M}athematical foundations and
  applications}, John Wiley \& Sons, Ltd., Chichester, 1990. \MR{1102677}

\bibitem{FHKL}
S.~Fischler, M.~Hussain, S.~Kristensen, and J.~Levesley, \emph{A converse to
  linear independence criteria, valid almost everywhere}, Ramanujan J.
  \textbf{38} (2015), no.~3, 513--528. \MR{3423011}

\bibitem{FKMS}
Lior Fishman, Dmitry Kleinbock, Keith Merrill, and David Simmons,
  \emph{Intrinsic {D}iophantine approximation on manifolds: general theory},
  Trans. Amer. Math. Soc. \textbf{370} (2018), no.~1, 577--599. \MR{3717990}

\bibitem{HussainKristensen}
M.~Hussain and S.~Kristensen, \emph{Badly approximable systems of linear forms
  in absolute value}, Unif. Distrib. Theory \textbf{8} (2013), no.~1, 7--15.
  \MR{3112202}

\bibitem{HussainKristensen2}
\bysame, \emph{Metrical results on systems of small linear forms}, Int. J.
  Number Theory \textbf{9} (2013), no.~3, 769--782. \MR{3043613}

\bibitem{Hussain}
Mumtaz Hussain, \emph{A note on badly approximable linear forms}, Bull. Aust.
  Math. Soc. \textbf{83} (2011), no.~2, 262--266. \MR{2784784}


\bibitem{HussainLevesley}
Mumtaz Hussain and Jason Levesley, \emph{The metrical theory of simultaneously
  small linear forms}, Funct. Approx. Comment. Math. \textbf{48} (2013),
  no.~part 2, 167--181. \MR{3100138}
  
    
  \bibitem{HSS}
Mumtaz Hussain, Johannes Schleischitz and David Simmons, \emph{The generalised Baker-Schmidt problem on hypersurface}, Pre-print: arXiv:1803.02314v1. To-appear: Int. Math. Res. Not. IMRN.

\bibitem{Jarnik3}
Vojt{\v e}ch Jarn\'ik, \emph{ {\"{U}ber die simultanen {D}iophantischen {A}pproximationen}}, Math. Z. \textbf{33} (1931), 505--543 (German). \MR{1545226}

\bibitem{Kleinbock4}
Dmitry Kleinbock, \emph{Badly approximable systems of affine forms}, J. Number
  Theory \textbf{79} (1999), no.~1, 83--102. \MR{1724255}

\bibitem{Levesley}
J.~Levesley, \emph{A general inhomogeneous {J}arnik-{B}esicovitch theorem}, J.
  Number Theory \textbf{71} (1998), no.~1, 65--80. \MR{1631030}

\bibitem{MHMK}
Seyyed~Hassan Mahboubi, Mumtaz Hussain, Abolfazl~S Motahari, and Amir~K
  Khandani, \emph{Layered interference alignment: Achieving the total {DOF} of
  {MIMO} {X}-channels}, Pre-print:arXiv:1412.7188 (2014).

\bibitem{McMullen}
Curt McMullen, \emph{Winning sets, quasiconformal maps and {D}iophantine
  approximation}, Geom. Funct. Anal. \textbf{20} (2010), no. 3, 726--740. \MR{2720230}

\bibitem{OrdUri}
Or~Ordentlich and Uri Erez, \emph{Precoded integer-forcing universally achieves
  the {MIMO} capacity to within a constant gap}, IEEE Trans. Inform. Theory
  \textbf{61} (2015), no.~1, 323--340. \MR{3299969}

\bibitem{Schmidt1}
Wolfgang~M. Schmidt, \emph{On badly approximable numbers and certain games},
  Trans. Amer. Math. Soc. \textbf{123} (1966), 27--50. \MR{0195595}

\bibitem{Schmidt2}
\bysame, \emph{Badly approximable systems of linear forms}, J. Number Theory
  \textbf{1} (1969), 139--154. \MR{0248090}

\bibitem{Sprindzuk}
Vladimir~G. Sprind{\v z}uk, \emph{Metricheskaya teoriya diofantovykh
  priblizheni\u\i ({M}etric theory of {D}iophantine approximations}, Izdat.
  ``Nauka'', Moscow, 1977 (Russian). \MR{0548467}

\bibitem{Sternberg}
Shlomo Sternberg, \emph{Lectures on differential geometry}, Prentice-Hall,
  Inc., Englewood Cliffs, N.J., 1964. \MR{0193578}

\end{thebibliography}

\end{document}